%% file: main.tex
\documentclass[a4paper]{article}
\usepackage[utf8]{inputenc}

\usepackage[pages=all, color=black, position={current page.south}, placement=bottom, scale=1, opacity=1, vshift=5mm]{background}
\SetBgContents{
}      %

\usepackage[margin=1in]{geometry} %
\usepackage{amsmath}
\usepackage{diagbox}
\usepackage{times}
\usepackage{subfig}
\usepackage{bm}
\usepackage{algorithm}
\usepackage{algpseudocode}
\usepackage{tabularx}
\usepackage{array}
\usepackage{multirow}
\usepackage{booktabs}
\usepackage{comment}
\usepackage{float}
\usepackage{stmaryrd}
\usepackage{authblk}
\usepackage{amsthm}
\usepackage{pgf-umlcd}
\usepackage{amssymb}
\usepackage{graphicx}
\usepackage{pgfplots}
\usepackage{todonotes}

\usepackage{placeins} %
\usepackage{pifont}%
\newcommand{\cmark}{\ding{51}}%
\newcommand{\xmark}{\ding{55}}%
\usepackage{siunitx}
\usepackage[export]{adjustbox}
\usepackage{pgfplotstable}
\usepackage{rotating}

\theoremstyle{definition}
\newtheorem{definition}{Definition}[section]

\DeclareMathOperator{\spn}{span}

\newcommand{\Pcf}{\mathcal{P}^{(c,f)}}
\usepgfplotslibrary{colorbrewer}
\pgfplotsset{
    cycle list/Dark2-7,
    cycle multiindex* list={
        mark list*\nextlist
        Dark2-7\nextlist
        },
}
\pgfplotsset{every axis plot/.append style={thick}}
\usepgfplotslibrary{groupplots,external}
\usetikzlibrary{er,positioning,external}
\tikzset{external/only named=true}
\usepackage{multicol}
\usepackage{mathtools}
\usepackage{caption} %

\usepackage{hyperref}
\hypersetup{
    colorlinks = true,
}

\usepackage[sort&compress,numbers,square]{natbib}
\DeclareMathOperator\supp{supp}
\usepackage{graphicx}
\graphicspath{{fig/}}
\usepackage{mathrsfs} %

\pgfplotsset{compat=1.9}

\DeclareMathOperator*{\argmin}{arg\,min}

\begin{document}
\title{Matrix-free implementation of the non-nested multigrid method}

\author[1,2]{Marco Feder}
\author[2]{Luca Heltai}
\author[3]{Martin Kronbichler}
\author[4]{Peter Munch}
\affil[1]{MathLab, Mathematics Area, SISSA}
\affil[2]{Numerical Analysis Group, University of Pisa}
\affil[3]{Numerical Mathematics, Ruhr University Bochum}
\affil[4]{Institute of Mathematics, Technical University of Berlin}

\maketitle

\begin{abstract}
\noindent Traditionally, the geometric multigrid method is used with nested levels. However, the construction of a suitable
hierarchy for very fine and unstructured grids is, in general, highly non-trivial. In this scenario, the non-nested multigrid method
could be exploited in order to handle the burden of hierarchy generation, allowing some flexibility on the
choice of the levels. We present a parallel, matrix-free, implementation of the non-nested
multigrid method for continuous Lagrange finite elements, where each level may consist of independently partitioned triangulations. Our algorithm has
been added to the multigrid framework of the C++ finite-element library \textsc{deal.II}.
Several 2D and 3D numerical experiments are presented, ranging from Poisson problems to linear elasticity. We test the robustness and performance of the proposed implementation with different polynomial degrees and geometries.
\end{abstract}

\section{Introduction}\label{sec:intro}

Scientific and industrial applications often require the solution of differential problems of the form
\begin{equation}
    \mathcal{A}u=f
\end{equation}
where $\mathcal{A}$ and $f$ are a given elliptic operator and source term, respectively, and $u$ the solution field of interest. Efficient
solvers relying on the finite element method (FEM) applied to complicated geometries often require large meshes with a large number of unknowns and adaptive mesh refinement in order 
to capture the behavior of the solution in some regions or to resolve the geometry of the domain. In order to solve the resulting large and sparse linear systems of equations,
fast and robust solvers and preconditioners have been developed.
Among various possibilities, multigrid
methods constitute one of the most efficient techniques for elliptic and parabolic problems~\cite{MGOverview}. Multigrid methods apply simple iterative schemes, which are effective at damping the high-frequency content of the iteration error and are called smoothers, on a hierarchy of resolution levels in order to reduce all error components quickly.

Multigrid methodologies come in various flavors, and the choice of variant strongly depends on the underlying mesh and finite-element space. Globally refined meshes
generated out of modestly-sized coarse grids are ideal candidates for the geometric multigrid method, which uses each mesh as a level
for the multigrid algorithm. In this context, the construction of coarse meshes out of a finer one is a simple task and can be exploited by element agglomeration when the
fine grid is structured, or even utilize the refinement information if a mesh is created by hierarchical refinement~\cite{CHKK,MHPSK}. On the other hand, for very fine and unstructured meshes on complicated geometries, stemming for instance from CAD (computer-assisted-design) models, it is
highly non-trivial to generate the suitable hierarchy of levels needed by the multigrid algorithm, representing a subject of active research.
In such cases, one could fall back to algebraic multigrid (AMG) \cite{StubenAMG} or to \emph{non-nested} multigrid variants \cite{BPX,Bittencourt}.
Multigrid methods for non-nested hierarchies have recently gained renewed attention in the context of polygonal-based methods such as the Virtual Element Method (VEM) and the Hybrid High Order Method (HHO), where the presence of
polyhedra induces non-nestedness between levels and hence an ad-hoc definition of the transfer operators has to be carried out, see e.g. ~\cite{HHOnonnestedMG,AntoniettiVEM} and references therein. For Discontinuous Galerking methods,
both nested (~\cite{AntoniettiNested}) and non-nested (~\cite{AntoniettiPennesi}) approaches are possible. Indeed, recent works propose agglomeration algorithms to generate polytopic hierarchies out of general meshes~\cite{feder2024r3mgrtreebasedagglomeration,antonietti2024polytopalmeshagglomerationgeometrical}, in order
to be used within multilevel methods. 

However, due to the lack of consolidated, automatic and high-quality coarsening strategies for complicated tetrahedral or hexahedral meshes, it is typical to build each level as an independent remeshing of the problem geometry by changing the mesh size inside a chosen mesh generator.

From the implementation standpoint, the main difficulty in the design of a geometric multigrid method for non-matching levels lies 
in the implementation of the intergrid operators. In \cite{KrauseDickopf}, different options have been thoroughly compared
both in terms of accuracy and computational cost for grids of moderate size. The necessary tasks for designing such transfer operators, as
for instance a spatial search procedure to identify which elements of the two meshes are intersecting, are shared by a large class of non-matching finite element techniques such as the fictitious domain method, the immersed boundary method,
X-FEM and Lagrange multiplier formulations \cite{GlowinskiPanPeriaux1994,peskin_2002} as well as particle-like methodologies \cite{golshan2021lethedem,JOACHIM2023112189}.
Most of these approaches require the computation of coupling operators representing the interaction of two unrelated meshes over which suitable finite element spaces and physics are defined, and the inaccurate computation of the entries of such a matrix
may in some cases hinder the convergence properties of the method \cite{HeinzMunchKaltenbacher,Laughton2021}. In \cite{KrauseZulian}, a parallel
technique for the variational transfer between unstructured meshes arbitrarily distributed among different processors is developed. 
Despite the different contexts and applications of the works mentioned so far, a shared feature is the explicit storage of the
coupling operator as a sparse and generally rectangular matrix that is later used in matrix-vector products.

In contrast, we propose a completely matrix-free implementation of a non-nested geometric multigrid solver and of its transfer operator. It is applicable for low and
moderately high polynomial degrees. For this purpose, we have extended
the multigrid infrastructure of \textsc{deal.II}~\cite{CHKK,MHPSK}, which has been
developed for nested meshes and optimized in the last decade~\cite{Kronbichler2018}. We compare the computational properties with geometric multigrid as well as with AMG and polynomial multigrid, indicating
good performance and robust behavior of the new implementation.

\subsection*{Main contribution}
\label{subsec:contrib}
To the best of the authors' knowledge, no matrix-free realizations of the multigrid method for finite-element discretizations for non-nested levels has been previously considered. In this work,
we will employ continuous Lagrangian finite elements and pointwise interpolation as transfer operator. 
The algorithms described in this work have been integrated into the widely used C++ finite-element library
\textsc{deal.II} \cite{dealIIdesign, dealII95} including
a tutorial\footnote{open pull request: \url{https://github.com/dealii/dealii/pull/17919}}; however, it should be noted that
the main building blocks are generic and could be inserted in other finite element frameworks which support the MPI standard \cite{GROPP1996789}, matrix free evaluation of weak forms on generic points in reference coordinates
as well as efficient geometric search capabilities. The current implementation supports simplex meshes and discontinuous nodal spaces, which are not considered by the analysis in this publication.

The rest of the work is organized as follows. In Section \ref{sec:nnmgmethod}, the non-nested multigrid method and its relevant
properties are recalled. Section \ref{sec:implementation} focuses on implementation details related to the matrix-free transfer operator between two arbitrarily
overlapping and distributed grids, and Section \ref{sec:numerical_experiments} presents a series of 2D and 3D numerical experiments that validate the described implementation and performance results. Eventually, some
examples from classical Finite Element Analysis are shown. Section~\ref{sec:performance} presents a performance analysis, discussing the main differences compared to classical nested multigrid algorithms. Conclusions are drawn in Section \ref{sec:conc}.

\input{chapters/non_nested_mg}

\input{chapters/implementation_details}

\input{chapters/sanity_checks}

\input{chapters/applications_fea}

\input{chapters/performance_evaluation}

\FloatBarrier
\section{Conclusions and outlook}\label{sec:conc}

\noindent We have presented a matrix-free and distributed-memory implementation of the non-nested multigrid method, which gives
flexibility for the levels that can be employed, allowing them to be arbitrarily overlapping and distributed among processors. We have shown its robustness
through an extensive set of examples by varying grids, equations, polynomial degrees, and the number
of levels. We confirmed its reliability also in the case of non-trivial three-dimensional geometries.
Building on top of the~\textsc{deal.II} finite element library, it has been possible to integrate our implementation
with highly-optimized and state-of-the-art matrix-free evaluation kernels. Finally, we have carried out
a breakdown of the different components of our pipeline, showing the computational cost associated
to the different phases and pointing to possible future improvements. Most of the infrastructure has been adapted in order to work seamlessly with simplex meshes, except
for exploiting the tensor product structure. Given the gain in flexibility with the non-nested multigrid infrastructure, more aggressive coarsening strategies beyond the classical 2:1 ratio are conceivable for future research. Similarly, extension of the transfer operators with algorithms beyond the simple polynomial-injection appears to be worthwhile.

\section*{Acknowledgments}

\noindent The authors acknowledge collaboration with Maximilian Bergbauer, Niklas Fehn, and Johannes Heinz as well as the deal.II community

\noindent LH acknowledges the MIUR Excellence Department Project awarded to the Department of Mathematics, University of Pisa, CUP I57G22000700001.

\noindent LH and MF acknowledge the partial support of the grant MUR PRIN 2022 No. 2022WKWZA8 “Immersed methods for multiscale and multiphysics problems (IMMEDIATE)”.

\noindent LH and MF are members of Gruppo Nazionale per il Calcolo Scientifico (GNCS) of Istituto Nazionale di Alta Matematica (INdAM).

\noindent MK and LH acknowledge partial support of the European Research Council (ERC) under the European Union's Horizon 2020 research and innovation programme (call HORIZON-EUROHPC-JU-2023-COE-03, grant agreement No. 101172493 ``dealii-X'').

\noindent MK and PM acknowledge partial support of the German Ministry of Education
and Research through the project ``PDExa: Optimized software methods for solving partial differential equations on exascale supercomputers'', grant agreement no. 16ME0637K, and the European Union -- NextGenerationEU.

\bibliography{refs}
\end{document}

%% file: chapters/non_nested_mg.tex
\section{Non-nested multigrid method}
\label{sec:nnmgmethod}
The main goal of this section is to describe the prolongation and restriction operators used to transfer between two consecutive grids. We refer the reader to \cite{HackbuschMG} for an extensive presentation of the
multigrid method and to \cite{BPX} for the convergence analysis of multigrid methods with non-nested spaces and non-inherited bilinear forms applied to elliptic problems.  For the sake of clarity, we consider the classical form of a variational problem: \emph{find} $u \in V$ \emph{such that} 

\begin{equation}
    \label{eqn:problem}
    a(u,v) = b(v) \qquad \forall v \in V,
\end{equation}
where $V$ is a suitable functional space for which bilinear and linear forms of the problem fulfill classical existence and uniqueness results defined on a domain $\Omega \subset \mathbb{R}^d$, $d=2,3$ with Lipschitz continuous boundary.
In view of the forthcoming application of multilevel techniques, let $\{ \mathcal{T}_{l}\}_{l=1,\ldots,L}$ denote a family of non-nested and shape-regular partitions for $\Omega$, each one characterized by disjoint open elements $K$ with diameter $h_K$, such that
$\overline{\Omega}=\bigcup_{K \in \mathcal{T}_l}\overline{K}$. The subscript $l \in \{1,\ldots, L\}$ is used to index multigrid levels and we indicate with $h_l = \max_{K \in \mathcal{T}_l} h_K$ the mesh size of the $l$-th level. In our notation, level $L$ consists of the finest mesh with size $h_L$, whereas $l=1$ represents the
coarsest level in the hierarchy. The coarser level $\mathcal{T}_{l}$ is independent of $\mathcal{T}_{l+1}$, with the only refinement constraint that there exists a constant $C>0$ independent of the discretization
parameters such that
\begin{equation}
    C h_{l} \leq h_{l+1} \leq h_{l} \qquad l \in \{1,\ldots,L-1\}.
\end{equation}
All grids consist of tensor product elements (quadrilaterals or hexahedra) and to each one of the levels we associate a set of degrees of freedom (DoFs) $\mathcal{D}^{l}$ and a conforming Lagrangian finite element space
$$V_l = \{v \in V: v_{|K} \in \mathcal{Q}^{p}(K), \forall K \in \mathcal{T}_l \} \subset V,$$ with $\mathcal{Q}^{p}(K)$ the space of polynomials of degree $p$ in each variable. In particular, the fact that $\mathcal{T}_{l} \not \subset \mathcal{T}_{l+1}$ induces the sequence of finite element spaces to be non-nested: $V_l \not \subset V_{l+1}$. The approximation of \eqref{eqn:problem} by finite elements and multigrid techniques requires solving several variational problems of the form
\begin{equation}
    \label{eqn:discrete_problem}
    a(u_l,v)=b(v), \qquad \forall v \in V_l.
\end{equation}
By introducing discrete operators $A_l: V_l \rightarrow V_l$ and $f_l:V_l \rightarrow \mathbb{R}$ defined by
\begin{equation}
    \bigl(A_l w,v\bigr)_{L^2(\Omega)} = a(w,v) \qquad \forall w,v \in V_l,
\end{equation}
and
\begin{equation}
    \bigl(f_l,v\bigr)_{L^2(\Omega)} = b(v) \qquad \forall v \in V_l,
\end{equation}
problem \eqref{eqn:discrete_problem} can be written as follows: \emph{find} $u_l \in V_l$ such that
\begin{equation}
    \label{eqn:linear_system}
    A_l u_l= f_l.
\end{equation}
The standard nodal interpolant for Lagrangian finite elements $I_l: C^0(\Omega) \rightarrow V_l$,
\begin{equation}
    \label{eqn:pointwise}
    u \mapsto I_l u \coloneqq \sum_{i \in \mathcal{D}^l} u(\bm{x}_i) \varphi_{l}^{i},
\end{equation}
can be employed to derive a transfer operator $\Pcf$ from the coarse to the fine grid for a given generic function $u \in C^{0}(\Omega)$, the space of continuous functions in $\Omega$. In~\eqref{eqn:pointwise}, $\varphi_{l}^{i}$ denotes the basis function associated to the $i$-th degree of freedom
in the mesh $\mathcal{T}_l$ located at the support point $\bm{x}_i$. Since we are employing Lagrangian nodal elements, each support point is associated to a unique degree of freedom. This does not hold true with (nodal) discontinuous Galerkin discretizations, where multiple degrees of freedom insist on
the same support point. We note that the point-to-point interpolation is a simple choice, and other strategies like mortar-based algorithms would also be possible~\cite{HeinzMunchKaltenbacher,Laughton2021}.
The residual associated to the approximate solution $\tilde{u}_l \in V_l$ is defined as
\begin{equation}
    \label{eqn:residual}
    r_l(v)\coloneqq b(v)-a(\tilde{u}_l,v) = \sum_{i \in \mathcal{D}^l} \bigl[ b(\varphi_{l}^{i}) - a(\tilde{u}_l,\varphi_{l}^{i}) \bigr] v(\bm{x}_i),
\end{equation}
for all $v \in V_l$, $l \in\{1,\ldots,L\}$. It is natural to associate a functional $\tilde{r}_l \in V_{l-1}^{'}$ defined by
\begin{align}
    \label{eqn:residual_dual}
    \tilde{r}_l: V_{l-1} \rightarrow \mathbb{R}, \qquad \tilde{r}_{l}(v_{l-1}) \coloneqq r_l(I_l v_{l-1}).
\end{align}
It is not difficult to verify that $\tilde{r}_l$ is a linear and bounded functional, the latter thanks to the coercivity assumption on $a$. For a given $v_{l-1} \in V_{l-1}$, based on employing \eqref{eqn:pointwise}, \eqref{eqn:residual}, and \eqref{eqn:residual_dual}, we have:
\begin{equation}
    \label{eqn:residual_interp}
    \begin{split}
        r_l(I_l v_{l-1}) &= b(I_l v_{l-1}) - a(\tilde{u}_l, I_l v_{l-1})  \\
        &= \sum_{i \in \mathcal{D}^l} \bigl[ b(\varphi_{l}^{i}) - a(\tilde{u}_l,\varphi_{l}^{i}) \bigr] v_{l-1}(\bm{x}_i)  \\
        &= \sum_{i \in \mathcal{D}^l} \Biggl( \bigl[ b(\varphi_{l}^{i}) - a(\tilde{u}_l,\varphi_{l}^{i}) \bigr] \sum_{j \in \mathcal{D}_{l-1}} v_{l-1}(\bm{x}_j) \varphi_{l-1}^{j}(\bm{x}_i) \Biggr) \\
        &= \sum_{j \in \mathcal{D}^{l-1}} \Biggl( \sum_{i \in \mathcal{D}_l} \bigl[ b(\varphi_{l}^{i}) - a(\tilde{u}_l,\varphi_{l}^{i}) \bigr]  \varphi_{l-1}^{j}(\bm{x}_i) \Biggr) v_{l-1}(\bm{x}_j) \\
        &= \sum_{j \in \mathcal{D}^{l-1}} \Biggl( \sum_{i \in \mathcal{D}_l} r_l(\varphi_l^i) \varphi_{l-1}^{j}(\bm{x}_i) \Biggr) v_{l-1}(\bm{x}_j).
    \end{split}
\end{equation}
This implies that the residual vector $r_l$ can be transferred from $\mathcal{T}_l$ to $\mathcal{T}_{l-1}$ with the matrix-vector product
\begin{equation}
    \label{eqn:residual2}
   r_{l-1} = \bigl( \Pcf\bigr)^T r_l,
\end{equation}
where 
\[
    \Pcf
    =
    \begin{bmatrix}
    \varphi_{l-1}^{1}(\bm{x}_1) & \varphi_{l-1}^{2}(\bm{x}_1) & \varphi_{l-1}^{3}(\bm{x}_1) & \dots  & \varphi_{l-1}^{|\mathcal{D}_{l-1}|}(\bm{x}_1) \\
    \varphi_{l-1}^{1}(\bm{x}_2) & \varphi_{l-1}^{2}(\bm{x}_2) & \varphi_{l-1}^{3}(\bm{x}_2) & \dots  & \varphi_{l-1}^{|\mathcal{D}_{l-1}|}(\bm{x}_2) \\
    \vdots & \vdots & \vdots & \ddots & \vdots \\
    \varphi_{l-1}^{1}(\bm{x}_{|\mathcal{D}_l|}) & \varphi_{l-1}^{2}(\bm{x}_{|\mathcal{D}_l|}) & \varphi_{l-1}^{3}(\bm{x}_{|\mathcal{D}_l|}) & \dots  & \varphi_{l-1}^{|\mathcal{D}_{l-1}|}(\bm{x}_{|\mathcal{D}_l|}) \\
\end{bmatrix}.\] 

\begin{figure}[!t]
    \centering
    \includegraphics{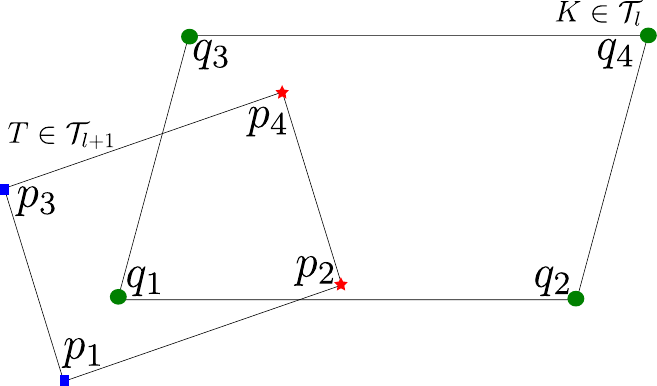}
    \caption{Two overlapping cells coming from consecutive levels. Green dots: DoFs associated to a $\mathcal{Q}^1$ element on the coarser cell $K$.
    Blue squares and red stars: DoFs associated to a $\mathcal{Q}^1$ element on the finer cell $T$. Red stars correspond to the DoFs $(\bm{p}_2,\bm{p}_4)$ falling inside $K \in \mathcal{T}_l$, while
    blue squares are the ones that are falling outside. Each $\bm{q}_i$ will evaluate on $(\bm{p}_2,\bm{p}_4)$ only. }
    \label{fig:overlapped_cells}
\end{figure}

\noindent Therefore, the generic entry of matrix $\Pcf$ corresponds to the evaluation of the $j$-th basis function associated with the coarser triangulation $\mathcal{T}_{l-1}$
on the $i$-th DoF attached to $\mathcal{T}_l$, and leads in general to a large and sparse rectangular matrix. 
If the point $\bm{x} \not \in \supp \varphi_{l-1}^j$, then the corresponding entry evaluates to $0$. Clearly, if parent-child relations between consecutive levels are not directly available, then ad-hoc search strategies must be performed in order to
identity the coarse-grid element $T \in \mathcal{T}_{l-1}$ s.t. $\bm{x}_i \in T$  for each finer point $\bm{x}_i$. This has the effect that the setup of the sparsity pattern of $\Pcf$ is more involved compared to the nested case, where only the knowledge of neighbors
is needed. A prototypical scenario is shown in Figure \ref{fig:overlapped_cells}, where some DoFs on the finer cells fall outside the coarser cell.
A more in-depth description of the information retrieval between two meshes and of the necessary data structures is given in Section \ref{sec:implementation}.
Given an initial guess $u_0 \in V_L$ and the numbers of pre-smoothing and post-smoothing steps
$(m_1,m_2) \in \mathbb{N}\times\mathbb{N}$, the non-nested multigrid V-cycle algorithm for the approximation of
$u_L$ is implemented in a recursive fashion, as summarized in Algorithm~\ref{alg:V_cycle}. Concrete choices for $\mathtt{PreSmoother}$, $\mathtt{CoarseGridSolver}$ and $\mathtt{PostSmoother}$
usually depend on the application at hand, discussed in Section~\ref{sec:numerical_experiments}. In the present work, we use the multigrid method as a preconditioner within
conjugate-gradient solvers, as it is often more robust~\cite{Trottenberg2001} than using it as a solver. As a critical remark, if the true geometry is \emph{complex}, it may not be possible to mesh it with a smaller number of cells while preserving
the relevant geometrical features.
\subsection*{Coarse grid solvers} %
Since the coarsest level may indeed not be really \emph{coarse}, higher order polynomials would already give a high number of degrees of freedom, affecting the choice of the $\mathtt{CoarseGridSolver}$ component in the multigrid pipeline. Indeed,
solving on the coarsest level could become non-negligible in terms of the total time to solution if an inappropriate coarse solver is used. In the case of complex meshes, using black-box AMG preconditioning might decay in efficiency for moderately high polynomial degree $p$~\cite{Kronbichler2018}, which is why we consider the usage of polynomial multigrid approaches~\cite{StadlerBirosHighOrderMG,brown2010efficient,OMALLEY2017177},
where the mesh size $h_l$ is kept constant on every level $l$ of the hierarchy, but the polynomial degree $p$ of the finite element space
is reduced incrementally from one level to the other. Among various possibilities, this work chooses to decrease the polynomial degree by one from one level to the next, i.e., the degree $p^{(c)}$ on the coarse mesh is defined as $p^{(c)} = p^{(f)}-1$, being $ p^{(f)}$ the degree defined on the finer mesh.

\begin{figure}[ht]
    \centering
    \includegraphics[width=6cm]{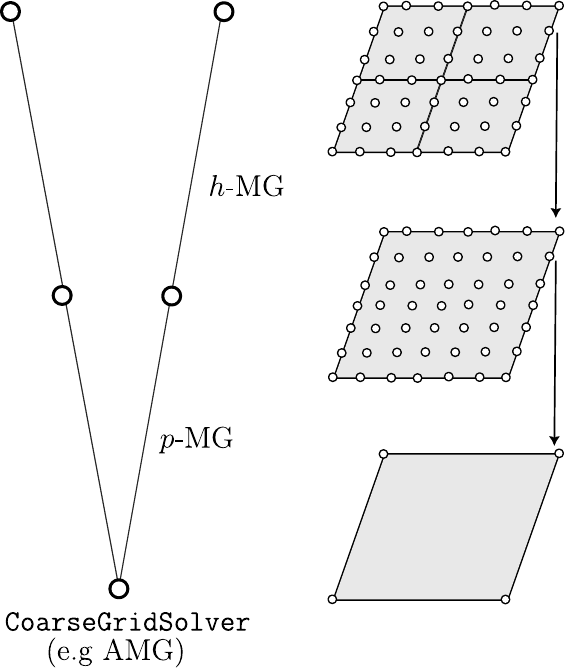}
    \hfill
    \includegraphics[width=6cm]{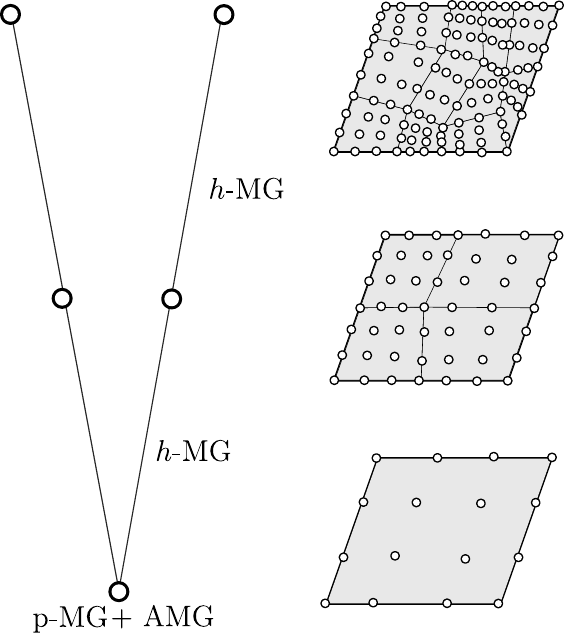}
    \caption{Schematic illustration of $hp$-multigrid scheme for a $\mathcal{Q}^3$ element. Support points corresponding to continuous Lagrangian elements are represented with white dots. 
    Left: Classical nested setting. Right: Non-nested $hp$ variant where the hierarchy of levels is \emph{non-matching}. Notice how the $\mathtt{CoarseGridSolver}$ is polynomial multigrid (p-MG).}
    \label{fig:hp_MG}
\end{figure}

The main benefit in our framework stems from the observation that whenever the coarsest level has many DoFs due to a high polynomial order $p$, we can change to
the conjugate-gradient method preconditioned by polynomial global coarsening. As we assume to work with linear elements at the coarsest level of the polynomial hierarchy, we can switch to AMG which is very competitive in that scenario~\cite{Kronbichler2018}.
Denoting with $\mathcal{T}_1$ the coarsest level over which we apply polynomial multigrid starting with degree $p$, then the sequence of finite element spaces $\{V_1^l\}_{l=1,\ldots,p}$
generated by the polynomial hierarchy is indeed \emph{nested} in the classical sense, as each space is defined on the same triangulation $\mathcal{T}_1$. In particular, highly optimized matrix-free kernels for polynomial transfer with nested levels are exploited \cite{Fehn2020,MHPSK}.
An illustration of this multigrid methodology and a comparison with the classical $hp$-multigrid is shown in Figure \ref{fig:hp_MG}. The outcome of this procedure can be regarded as an $hp$-multigrid scheme where the "$h$-" component may be non-nested. 

\begin{algorithm}[!t]
    \caption{One iteration of non-nested multigrid V-cycle on level $l \geq 2$ to solve $A_{l}x =f_{l}$}
    \begin{algorithmic}
        \State $x_0$ initial guess
        \If{$l=1$}
            \State $\delta_1\ \gets \mathtt{CoarseGridSolver}(A_{1},f_{1})$
        \Else
            \State $x_{l} \gets \mathtt{PreSmoother}\bigl(A_{l},0,f_{l},m_1\bigr)$
            \State $r_{l} \gets f_{l} - A_{l} x_{l} $
            \State $r_{l-1} \gets \mathtt{Restrictor}\bigl(r_{l}\bigr)$
            \State $\delta_{l-1} \gets \mathtt{V-cycle}\bigl(A_{l-1},r_{l-1},0,m_1,m_2\bigr)$
            \State $x_{l} \gets x_{l} + \mathtt{Prolongator}\bigl(\delta_{l-1}\bigr)$
            \State $x_{l} \gets \mathtt{PostSmoother}\bigl(A_{l},0,f_{l},m_2\bigr)$
        \EndIf
    \end{algorithmic}
    \label{alg:V_cycle}
\end{algorithm}

%% file: chapters/implementation_details.tex
\section{Implementation details}
\label{sec:implementation}
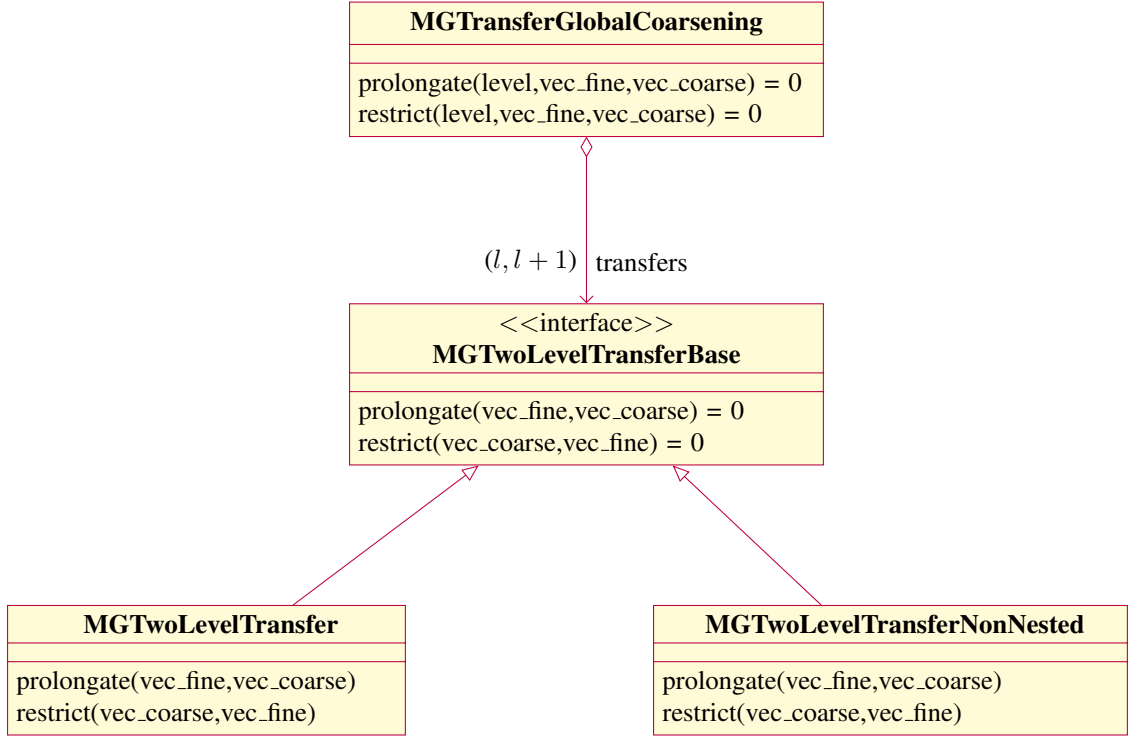
\begin{figure}[!t]
    \begin{tikzpicture}
        \begin{class}[text width=6cm]{MGTransferGlobalCoarsening}{0 ,0}
            \operation{prolongate(level,vec\_fine,vec\_coarse) = 0}
            \operation{restrict(level,vec\_fine,vec\_coarse) = 0}
        \end{class}
        
        \begin{interface}[text width=6.cm]{MGTwoLevelTransferBase}{0 ,-4}
            \operation{prolongate(vec\_fine,vec\_coarse) = 0}
            \operation{restrict(vec\_coarse,vec\_fine) = 0}
        \end{interface}
        \aggregation{MGTransferGlobalCoarsening}{transfers}{($l,l+1$)}{MGTwoLevelTransferBase}
        
        \begin{class}[text width=5.cm]{ MGTwoLevelTransfer
            }{ -5 , -8}
            \inherit {MGTwoLevelTransferBase}
            \operation{prolongate(vec\_fine,vec\_coarse)}
            \operation{restrict(vec\_coarse,vec\_fine)}
        \end{class}
        
        \begin{class}[text width=6.cm]{ MGTwoLevelTransferNonNested }{4 , -8}
            \inherit {MGTwoLevelTransferBase}
            \operation{prolongate(vec\_fine,vec\_coarse)}
            \operation{restrict(vec\_coarse,vec\_fine)}
        \end{class}
    \end{tikzpicture}
    \caption{UML diagram of transfer operators available in the $\mathtt{MGTransferGlobalCoarsening}$ framework in \textsc{deal.II}.
    The new abstract class delegates the implementation of the intergrid transfers to the derived class $\mathtt{MGTwoLevelTransfer}$
    (used in case of nested meshes) or to the new $\mathtt{MGTwoLevelTransferNonNested}$ in case of a non-nested multigrid method. Each two-level
     transfer object is specific to consecutive levels $l$ and $l+1$.}
    \label{fig:UML_diagram}
\end{figure}

Since each level $l$ covers the whole computational domain as in geometric global coarsening algorithms \cite{BeckerBraack,KANSCHAT20042437,BastianMG},
it is natural to implement the present framework in the already available global coarsening infrastructure present in the \textsc{deal.II} finite element library. The global coarsening infrastructure ($\mathtt{MGTransferGlobalCoarsening}$), extensively described in \cite{MHPSK} for nested grids, delegates the
actual transfer operations to a two-level implementation class ($\mathtt{MGTwoLevelTransfer}$). In order to allow a non-nested transfer within the same
interface, a new abstract base class ($\mathtt{MGTwoLevelTransferBase}$) has been introduced.
A UML diagram is displayed in Figure \ref{fig:UML_diagram}.

\begin{figure}[ht]
    \centering
    \subfloat[\centering Ball inside the cube]{{\includegraphics[width=7cm]{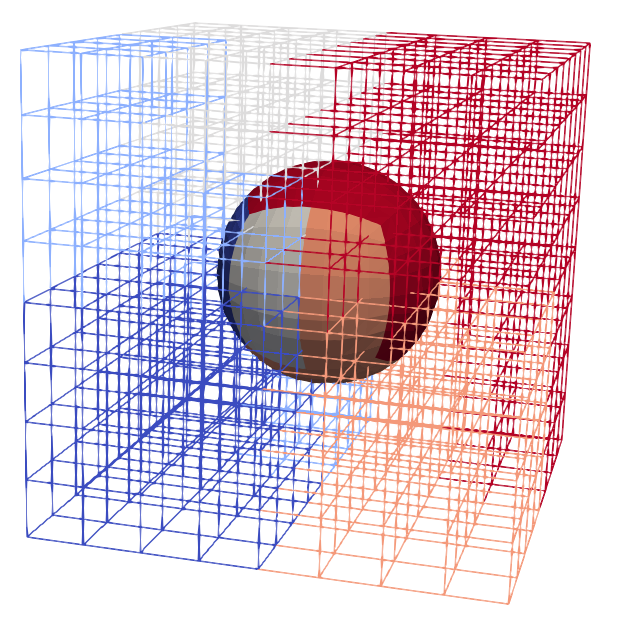} }}
    \subfloat[\centering Clip view]{{\includegraphics[width=7cm]{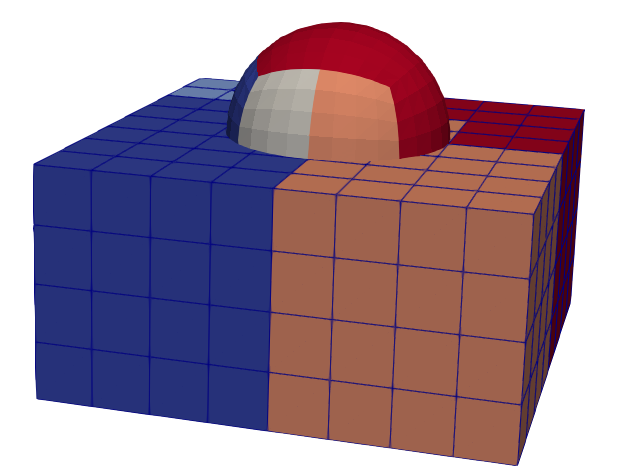} }}%
    \caption{Coupling between processes for two overlapped and distributed triangulations (each color represents a different rank). (a) The two partitioned triangulations. The cube is displayed with
    a wireframe view in order to show the inner ball. (b) Clip view to highlight the geometric intersection of mesh elements belonging to different processes. Notice that grids are discretizing two different geometries only for the purpose of showing the issue. In practice,
    they will both discretize the same domain $\Omega$.}
    \label{fig:distributed_trias}
\end{figure}

Among the several components required by a multigrid framework, there are two main ingredients that are detailed
in this section. First, due to the overlapping nature of the two grids, a \emph{geometric search} must be carried out to obtain the
mutual relationships between two consecutive levels $l$ and $l+1$, i.e., which cells on level $l$ own the DoFs of level $l+1$ and the associated reference positions.
In a distributed setting, where all triangulations are partitioned independently
(see Figure~\ref{fig:distributed_trias}), the relevant owner processes need to be determined by appropriate communication patterns.
Second, a matrix-free transfer operator has to evaluate the solution field
on arbitrarily located points
in the reference cell $[0,1]^d$ using the underlying polynomial basis. The same setting is well-established in distributed implementations of immersed methods, and we refer to \cite{KrauseZulian} for an extensive description of this task and details concerning load balancing.

\subsection{Distributed geometric search}
\label{subsubsec:Search procedure}
The evaluation of a discrete finite element field on a given list of physical points requires to determine the cell $K$ owning every point $\bm{x}$ along with its reference coordinates $\hat{\bm{x}}$. As described before, the major difficulty 
in the distributed case is that $K$ and $\bm{x}$ might belong to different processes. Consider a communicator consisting of $p$ processes, a list $\{\bm{x}_i\}_{i=1}^{N_q}$ of $N_q$ points locally owned by process $q \ne p$, and the task to \emph{find cells and corresponding ranks owning each $\bm{x}_i$}. A first coarse search is used to determine, in a cheap way, the candidate ranks owning each entity. A possible option uses a global communication step where
each process shares a rough description of the local portion of the domain it owns with other participants, allowing the global mesh description to be available to each process. Using the MPI standard, this can be achieved
via an $\mathtt{MPI\_Allgather}$ with vectors of Axis-Aligned Bounding Boxes (AABB) local to each rank. Thanks to a local tree data structure, for every element of $\{\bm{x}_i\}_{i=1}^{N_q}$ one can find possible owning cells and ranks. Once candidates have been
determined, a finer search is carried out by each process, checking whether it actually owns the points or not. The sequence of requests and answers can be realized efficiently by using consensus-based algorithms for dynamic sparse communications \cite{NBX}.
In order to avoid the quadratic complexity of $\mathtt{MPI\_Allgather}$ in the number of ranks, we optionally use a distributed tree search provided by the \textsc{ArborX} library \cite{ArborX} during the coarse search: each process builds the local tree
out of a local vector of AABB. Then, roots of all the local trees are used to create a second tree used for querying possible ownership of entities. Alternatively, other techniques that allow to determine ownership of possibly remote points based on forest-of-trees approaches are possible \cite{Burstedde_ParallelTree}. Search strategies of this kind are general and can be performed with other \emph{geometric entities} for which suitable simple predicates can be queried.
Once all the owners have been determined, reference positions $\{\bm{\hat{x}}_i\}_{i=1}^{N_q}$ associated to each $\{\bm{x}_i\}_{i=1}^{N_q}$ are computed by inversion of the transformation map from real to unit cell. For affine mappings, the inversion can be done explicitly, otherwise Newton-like methods are used, using tensor-product evaluation of a polynomial representation of the geometry, similar to the evaluation considered in the next subsection, see also~\cite{Bergbauer2024}.

\subsection{Efficient evaluation at reference positions}
The cell $K$ and reference position $\bm{\hat{x}}$ of each point are determined in a setup phase of the multigrid transfer operator. For the actual transfer operation, the remaining work is to evaluate the finite element solution at each point:

\begin{equation}
    \label{eqn:local_evaluation}
    u_K(\bm{x}) = \sum_{ 1 \leq i \leq N_K} {\varphi}_{i}\left(\hat{\bm{x}}_K\right) u_{K,i},
\end{equation} where $N_K$ is the number of degrees of freedom on cell $K$ and $\varphi_{i}$ the basis functions in reference coordinates with $\hat{\bm x}_K$ the reference position in cell $K$ for the point $\bm x$. For tensor-product shape functions in $\mathbb{R}^d$, the evaluation of a basis
function $\varphi$ at point $\hat{\bm{x}}$ can be written as follows:
\begin{equation}\label{eqn:tensor_product}
    {\varphi}(\hat{\bm{x}}) = \prod_{j = 1, \ldots, d} \varphi^{\text{1D}}_{j}\left(\hat{\bm{x}}^j\right),
\end{equation}
with $\varphi^{\text{1D}}_j(\cdot)$ being the one-dimensional shape function in direction $j$ and $\hat{\bm{x}}^j$ the $j$-th coordinate of point $\hat{\bm{x}}$. Combining Equations~\eqref{eqn:local_evaluation} and~\eqref{eqn:tensor_product}, we get
\begin{equation}
    \label{eqn:tensor_product_evaluation}
    u_K(\bm{x}) = \sum_{ 1 \leq k \leq N_{\text{DoFs}}^{\text{1D}}} \varphi^{\text{1D}}_k(\hat{\bm{x}}^3) \sum_{ 1 \leq j \leq N_{\text{DoFs}}^{\text{1D}}} \varphi^{\text{1D}}_j(\hat{\bm{x}}^2) \sum_{ 1 \leq i \leq N_{\text{DoFs}}^{\text{1D}}} \varphi^{\text{1D}}_i(\hat{\bm{x}}^1) u_{K,{ijk}},
\end{equation} using multi-indices $(i,j,k)$ that are related to the index $i$ in ~\eqref{eqn:local_evaluation} in a bijective way. Following the algorithm presented in~\cite{KRONBICHLER2012135} for a single point with optimizations from~\cite{Bergbauer2024}, the one-dimensional basis functions at $\hat{\bm x}^j$, $j=1,2,3$, are tabulated. With the tabulated arrays $N_{\bm{x}}^{1}, N_{\bm{x}}^{2}, N_{\bm{x}}^{3}$, one per direction with size $N_{\text{DoFs}}^{\text{1D}}$ each, equation ~\eqref{eqn:tensor_product_evaluation} can be written in tensor-product notation as
\begin{equation}
    u_K(\bm{x}) = N_{\bm{z}}^{3} \bigl( I_z \otimes N_{\bm{x}}^{2} \bigr) \bigl( I_z \otimes I_y \otimes N_{\bm{x}}^1 \bigr)u_K.
\end{equation}
The resulting computational complexity per point is $\mathcal{O}\left(\left(N_{\text{DoFs}}^{\text{1D}}\right)^d\right)$. 
Note that if the evaluation points within a cell have a tensor-product structure, classical sum-factorization techniques \cite{Orszag1980,KRONBICHLER2012135} combining evaluation steps for all points on a cell yield a final complexity of $\mathcal{O}\left(dN_{\text{DoFs}}^{\text{1D}}\right)$ per degree of freedom. For nested transfers, we use this strategy.

To summarize the computational properties, the present non-nested algorithm has a similar $\mathcal{O}\left(\left(N_{\text{DoFs}}^{\text{1D}}\right)^d\right)$ complexity as
matrix-based variants for the transfer. The crucial difference to previous methods is the fact that the identified computational complexity of unstructured matrix-free evaluation proposed in the present contribution happens on cached data. As a result, the memory access is reduced and the arithmetic intensity increases. This suggests an advantage of around one order of magnitude for contemporary hardware for higher polynomial degrees~\cite{Kronbichler2018,Bergbauer2024}.

\begin{figure}[ht]
    \centering
    \hfill
    \subfloat[\centering Gauss-Lobatto points for degree $p=4$.]{{\includegraphics{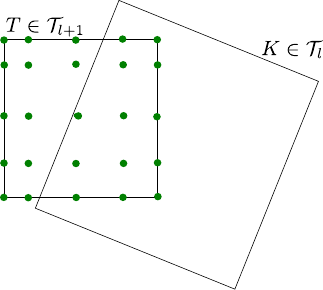} }}
    \hfill
    \subfloat[\centering Finer Gauss-Lobatto points seen from coarser cell $K$.]{{\includegraphics{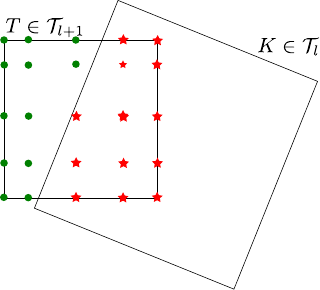} }}%
    \hfill
    \phantom{}
    \caption{(a) Gauss-Lobatto points for a quadrature rule of order $p=4$ on a cell $T \in \mathcal{T}_{l+1}$. (b) Evaluation points seen from the coarser cell $K \in \mathcal{T}_l$ (red stars) do not have a tensor-product structure.}%
    \label{fig:GaussLobatto}
\end{figure}

\subsection{Domains with curved boundaries}

\begin{figure}[!t]
  \centering
  \hfill
  \subfloat[\centering Non-conformity at the boundaries of different levels.]{{\includegraphics[scale=0.65]{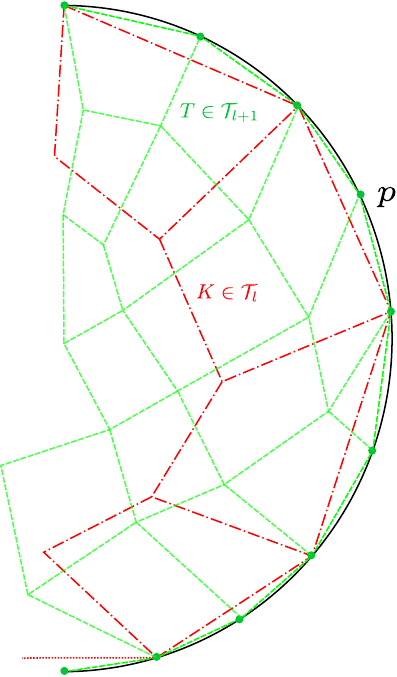} }}
  \hfill
  \subfloat[\centering Position of $F_K^{-1}(\bm{p})$ with respect to reference cell $\hat{K}$.]{{\includegraphics[width=5.5cm]{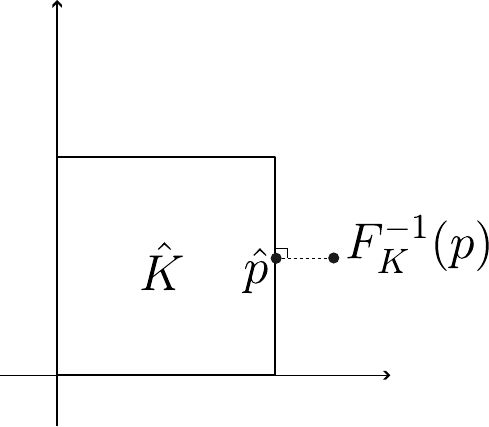} }}%
  \hfill
  \phantom{}
  \caption{Illustration of the situation occurring when the hierarchy discretizes a two-dimensional domain $\Omega$ with curved boundaries. (a) Red dash-dotted lines: elements of a coarser triangulation $\mathcal{T}_{l}$. Green dashed lines: boundary for a finer triangulation $\mathcal{T}_{l+1}$. Black solid line: exact representation of the boundary $\partial \Omega$. (b) Pre-image of point $\bm p$ through $F_K$ and its nearest point projection $\hat{\bm p} \in \overline{\hat{K}}$ defined according to \eqref{eqn:closest_point}.}%
  \label{fig:closest_projection}
\end{figure}

For simple geometrical shapes such as shells and cylinders, the placement of new vertices upon mesh refinement can be done through polar or cylindrical coordinates. In practice, more complex and realistic
shapes are modeled by means of CAD tools, which are internally exploiting B-splines and NURBS to represent most surfaces of interest with high accuracy. In consequence, the external generation
by meshing software of a sequence of refined grids that discretize an input complex geometry $\Omega$ gives rise to non-conformity on the boundaries of consecutive levels if $\partial \Omega$ is curved. Furthermore, this approach uses the geometrical information only at the pre-processing stage of meshing, not at later stages of the pipeline where it may be required to use high order representations of mappings from
real to unit cells. Indeed, even if each level (mesh) has boundary vertices \emph{placed correctly} on $\partial \Omega$, support points defining the shape functions on each boundary face will not lie on the manifold $\partial \Omega$.
The non-conformity implies that some entries of the matrix $\Pcf$, which should be non-zero, will vanish on points on the finer grid if not contained in the mesh involving coarser shape functions, resulting in an inaccurate intergrid transfer. For Dirichlet boundary
conditions on a portion $\Gamma^D \subset \partial \Omega$, this does not constitute a problem as DoFs on $\Gamma^D$ are already constrained. Conversely, for Neumann boundaries $\Gamma^N$ it is necessary 
to include all unconstrained DoFs on $\Gamma^N$ in order to provide an accurate transfer operator. This situation is illustrated in Figure \ref{fig:closest_projection}. Point $\bm p$ lies on a boundary face of $\mathcal{T}_{l+1}$
and falls outside the coarser level $\mathcal{T}_l$. During the geometric search procedure, by expanding local bounding boxes by a large enough tolerance, $\bm p$ is associated to the red element $K \in \mathcal{T}_l$. Denoting by $\bm F_K \colon \hat{K}\rightarrow K$ the unit-to-real map to cell $K$,
it holds that $\bm F_K^{-1}(\bm p) \not \in \overline{{\hat{K}}}$. Hence, shape functions are evaluated at point $\hat{\bm p} \in \overline{\hat{K}}$ nearest to $\bm F_K^{-1}(\bm p)$:
\begin{equation}
    \label{eqn:closest_point}
    \hat{\bm p} \coloneqq \argmin_{\hat{\bm x} \in \overline{\hat{K}}} d \left(\hat{\bm x},\bm F_K^{-1}(\bm p)\right),
\end{equation} with $d(\cdot,\cdot)$ the Euclidean distance on $\mathbb{R}^d$. Another approach, more robust and requiring fewer heuristics, is to propagate the geometrical information stored inside the CAD into the finite element computations by using state-of-the-art libraries such as \textsc{opencascade}~\cite{opencascade}
to interrogate the CAD model during the refinement process and the distribution of evaluation points on boundary faces \cite{geometry_paper}.

%% file: chapters/sanity_checks.tex
\section{Numerical experiments}
\label{sec:numerical_experiments}
This section presents several numerical experiments to verify the robustness of our algorithm with respect to 2D and 3D geometries and polynomial degrees by means of 2D and 3D Poisson problems with homogeneous Dirichlet boundary conditions and constant right-hand side $f=1$. Some sanity checks aimed
to check basic properties of the algorithm in well-known simple cases are shown first before considering actual non-matching levels. Continuous Lagrangian finite elements defined on quadrilaterals or hexahedra with polynomial degree $p$ ranging from $1$ to $4$ are employed. As quadrature formula, the Gauss-Legendre rule with $(p+1)^d$ points per cell is considered, with $d=2,3$. 
As discussed in Section~\ref{sec:implementation}, each level consists of a sequence of (possibly) distributed and unstructured triangulations meshed independently.
In 2D tests, unstructured levels are meshed by using the mesh generator \textsc{gmsh} \cite{Gmsh}. Concerning 3D tests,
the commercial software \textsc{Coreform Cubit} \cite{Cubit} has been adopted in order to generate high quality unstructured hexahedral meshes.
Finally, we conclude the section to applications of the proposed methodology with non-trivial geometries stemming from classical Finite Element Analysis (FEA): a complicated geometry is first read
from a CAD model and then meshed by an external software. We remark that in practice, a CAD file must first be repaired in order to be meshed. Removing the so-called small features from a CAD file
is necessary to mesh the model at hand and in general is a non-trivial operation that could require some time. After that, specific mesh-related parameters have to be carefully tuned, depending on the physical situation at hand. The grids generated
through this procedure are inherently \emph{non-nested}. The iterative solver is configured in the following way, as in~\cite{MHPSK}, with the exception that in the present work
a \emph{non-nested} multigrid preconditioner is employed in place of a nested one:

\begin{itemize}
    \item The conjugate-gradient solver is run until a reduction of the $l_2$-norm of the unpreconditioned residual by $10^4$ is reached. Such a loose tolerance is typical in many multigrid applications, for instance in time-dependent problems in computational fluid dynamics, for which good initial guesses are usually available by employing extrapolation or projection.
    \item The conjugate-gradient solver is preconditioned by a single V-cycle of the \emph{non-nested} multigrid method.
    \item Operations in the V-cycle are run with single-precision floating-point numbers, while conjugate-gradient is run in double precision, which increases the computational throughput while maintaining acceptable accuracy~\cite{KRONBICHLER2012135,Kronbichler2019}.
    \item A point Jacobi smoother employed within a Chebyshev iteration of degree 3 is used on every level, using eigenvalue estimates computed with 12 iterations of the Lanczos iteration.
    \item Two V-cycles of the \textsc{Trilinos ML} implementation of AMG \cite{TrilinosML} are used, in double precision, as coarse-grid solver. The used parameters are the same as the ones in Appendix C of \cite{MHPSK}.
\end{itemize}
In every test, we report the number of degrees of freedom attached to each level for different polynomial degrees. All the experiments have been performed on the Galileo100 Italian supercomputer\footnote{Using the deal.II master e55124254b with the Intel 2021.10.0.20230609 compiler and \textsc{-O3} and \text{-march=native} flags.}. Its compute nodes have two sockets (each one with 24 cores of Intel CascadeLake) and AVX-512 ISA extension, which allows 8 doubles or 16 floats to be processed per instruction.\footnote{For the technical specifications, see \href{https://www.hpc.cineca.it/systems/hardware/galileo100/}{https://www.hpc.cineca.it/systems/hardware/galileo100/}, retrieved on December 3, 2024.}

The wall-clock times (measured in seconds) shown in the upcoming tables consist of the actual solve time needed by the conjugate-gradient iteration. All results shown in this publication, along with the necessary meshes and CAD files, are available on a maintained GitHub repository \href{https://github.com/peterrum/dealii-multigrid}{https://github.com/peterrum/dealii-multigrid}.

\subsection{Application to nested meshes}

The present test is only meant to assess the consistency of the algorithm. A sequence of structured and 
globally refined nested meshes discretizing $[-1,1]^2$ is considered. Since $V_{l-1} \subset V_l$, the transfer operator
$\Pcf$ coincides with the classical injection from $V_{l-1}$ to $V_l$. As a matter
of fact, the numerical results obtained both with the nested (i.e. global coarsening) and the non-nested method are identical, as reported in
Table \ref{tab:nested_2D}. The same sanity check is repeated for a sequence of nested levels discretizing the unit cube $[-1,1]^3$. Results are reported in Table \ref{tab:nested_3D}. 

\begingroup
\renewcommand{\arraystretch}{1.1}
    \begin{table}[!t]
        \centering
        \resizebox{\columnwidth}{!}{%
        \begin{tabular}{l|cccr|cccr|cccr|cccr}
            \toprule
            \emph{l}& \multicolumn{4}{|c}{$p=1$} & \multicolumn{4}{|c}{$p=2$} & \multicolumn{4}{|c}{$p=3$} & \multicolumn{4}{|c}{$p=4$}  \\
            \cmidrule(lr){2-5}\cmidrule(lr){6-9}\cmidrule(lr){10-13}\cmidrule(lr){14-17}
            &$\text{\#i}_{\text{NN}}$ &$\text{\#i}_{\text{GC}}$ & $\text{\#i}_{\text{AMG}}$ & \#DoFs &$\text{\#i}_{\text{NN}}$ &$\text{\#i}_{\text{GC}}$  & $\text{\#i}_{\text{AMG}}$ & \#DoFs &$\text{\#i}_{\text{NN}}$ &$\text{\#i}_{\text{GC}}$ & $\text{\#i}_{\text{AMG}}$ & DoFs &$\text{\#i}_{\text{NN}}$ &$\text{\#i}_{\text{GC}}$ & $\text{\#i}_{\text{AMG}}$ & \#DoFs \\
            \midrule
            2 & 3 & 3 & 1 & 9   & 3 & 3 & 3 & 25    & 3 & 3 & 4  & 49    & 3 & 3 & 6  & 81  \\
            3 & 3 & 3 & 3 & 25  & 3 & 3 & 5 & 81    & 3 & 3 & 8  & 169   & 3 & 3 & 12 & 289  \\
            4 & 3 & 3 & 4 & 81  & 3 & 3 & 7 & 289   & 3 & 3 & 10 & 625   & 3 & 3 & 14 & 1 089  \\
            5 & 3 & 3 & 6 & 289 & 3 & 3 & 8 & 1 089 & 3 & 3 & 14 & 2 401 & 3 & 3 & 20 & 4 225  \\
            \bottomrule
        \end{tabular}
        }
        \caption{Number of iterations for the \emph{nested} sanity check in two dimensions. Legend: $\text{\#i}_{\text{NN}}$: number of iterations for non-nested multigrid. $\text{\#i}_{\text{GC}}$: number of iterations for the global coarsening algorithm. $\text{\#i}_{\text{AMG}}$: number of iterations required by AMG.}
        \label{tab:nested_2D}
    \end{table}
\endgroup 

\begingroup
\renewcommand{\arraystretch}{1.1}
    \begin{table}[!t]
        \centering
        \resizebox{\columnwidth}{!}{%
        \begin{tabular}{l|cccr|cccr|cccr|cccr}
            \toprule
            \emph{l}& \multicolumn{4}{|c}{$p=1$} & \multicolumn{4}{|c}{$p=2$} & \multicolumn{4}{|c}{$p=3$} & \multicolumn{4}{|c}{$p=4$}  \\
            \cmidrule(lr){2-5}\cmidrule(lr){6-9}\cmidrule(lr){10-13}\cmidrule(lr){14-17}
            &$\text{\#i}_{\text{NN}}$ &$\text{\#i}_{\text{GC}}$ & $\text{\#i}_{\text{AMG}}$ & \#DoFs &$\text{\#i}_{\text{NN}}$ &$\text{\#i}_{\text{GC}}$  & $\text{\#i}_{\text{AMG}}$ & \#DoFs &$\text{\#i}_{\text{NN}}$ &$\text{\#i}_{\text{GC}}$ & $\text{\#i}_{\text{AMG}}$ & DoFs &$\text{\#i}_{\text{NN}}$ &$\text{\#i}_{\text{GC}}$ & $\text{\#i}_{\text{AMG}}$ & \#DoFs \\
            \midrule
            2 & 3 & 3 & 1 & 27    & 3 & 3 & 3  & 125    & 3 & 3 & 10 & 343     & 3 & 3 & 20  & 729       \\
            3 & 3 & 3 & 4 & 125   & 3 & 3 & 10 & 729    & 3 & 3 & 20 & 2 197   & 3 & 3 & 33  & 4 913     \\
            4 & 3 & 3 & 5 & 729   & 3 & 3 & 12 & 4 913  & 3 & 3 & 31 & 15 625  & 3 & 3 & 47  & 35 937    \\
            5 & 3 & 3 & 6 & 4 913 & 3 & 3 & 16 & 35 937 & 3 & 3 & 73 & 117 649 & 3 & 3 & 100 & 274 625    \\
            \bottomrule
        \end{tabular}
        }
        \caption{Number of iterations for the \emph{nested} sanity check in three dimensions. Legend: $\text{\#i}_{\text{NN}}$: number of iterations for non-nested multigrid. $\text{\#i}_{\text{GC}}$: number of iterations for the global coarsening algorithm. $\text{\#i}_{\text{AMG}}$: number of iterations required by AMG.}
        \label{tab:nested_3D}
    \end{table}
\endgroup

\subsection{L-shaped domains} \label{subsub:L_shaped}
\subsubsection{2D unstructured L-shaped domains} \label{subsubsec:sanity_check_2D_l_shaped}
\begin{figure}[!t]
    \subfloat{\includegraphics[width=0.21\linewidth]{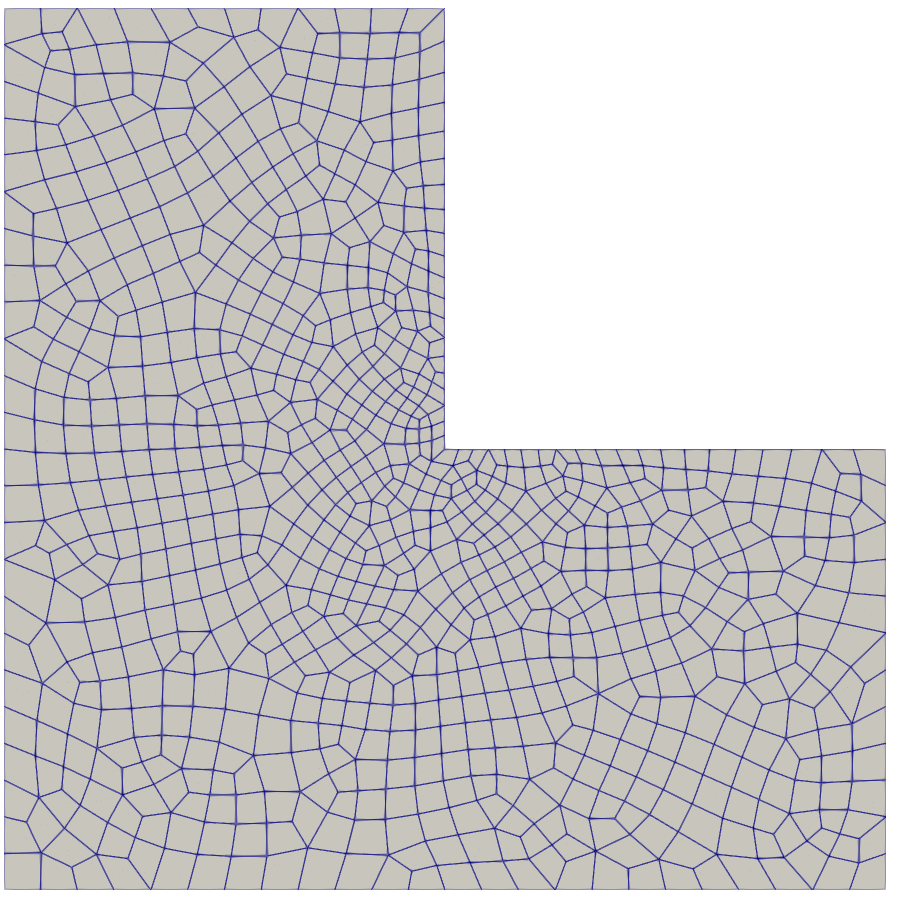}}\qquad%
    \subfloat{\includegraphics[width=0.21\linewidth]{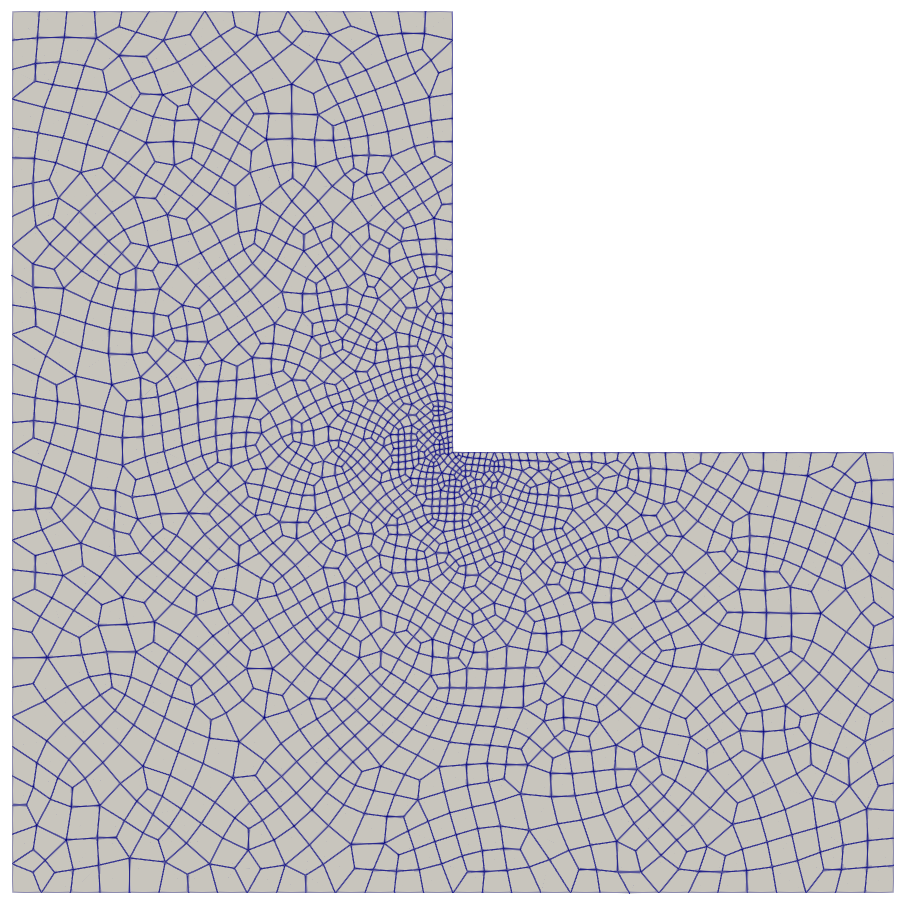}}\qquad%
    \subfloat{\includegraphics[width=0.21\textwidth]{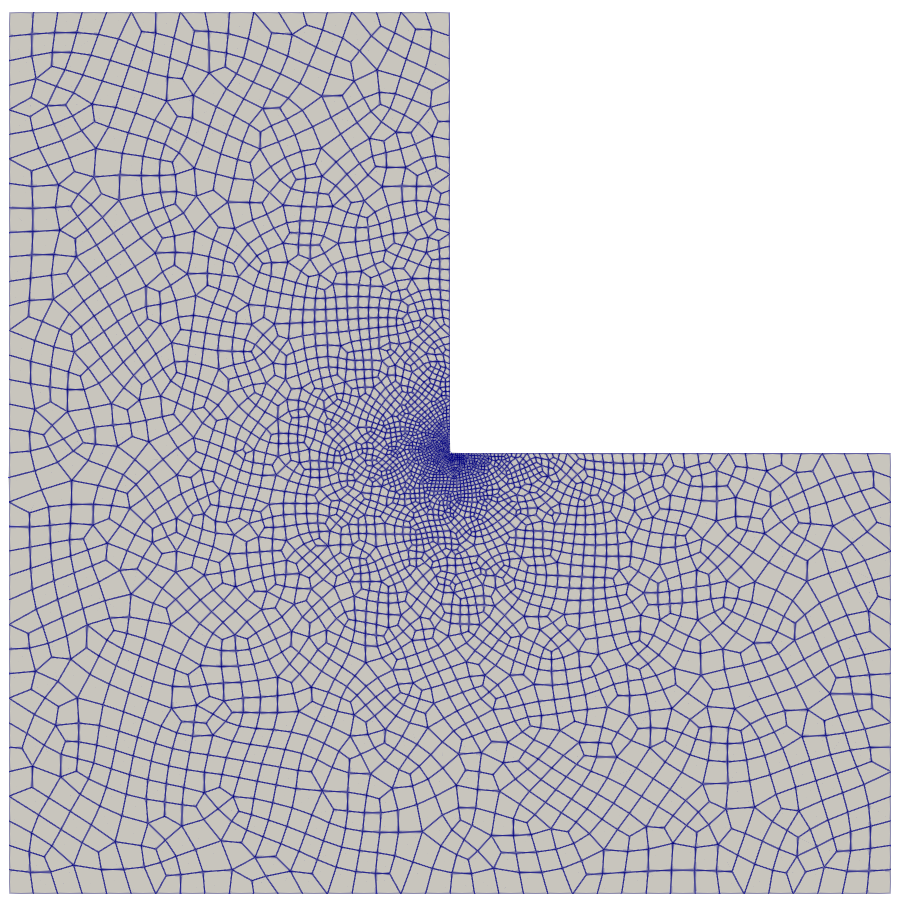}}\qquad%
    \subfloat{\includegraphics[width=0.21\textwidth]{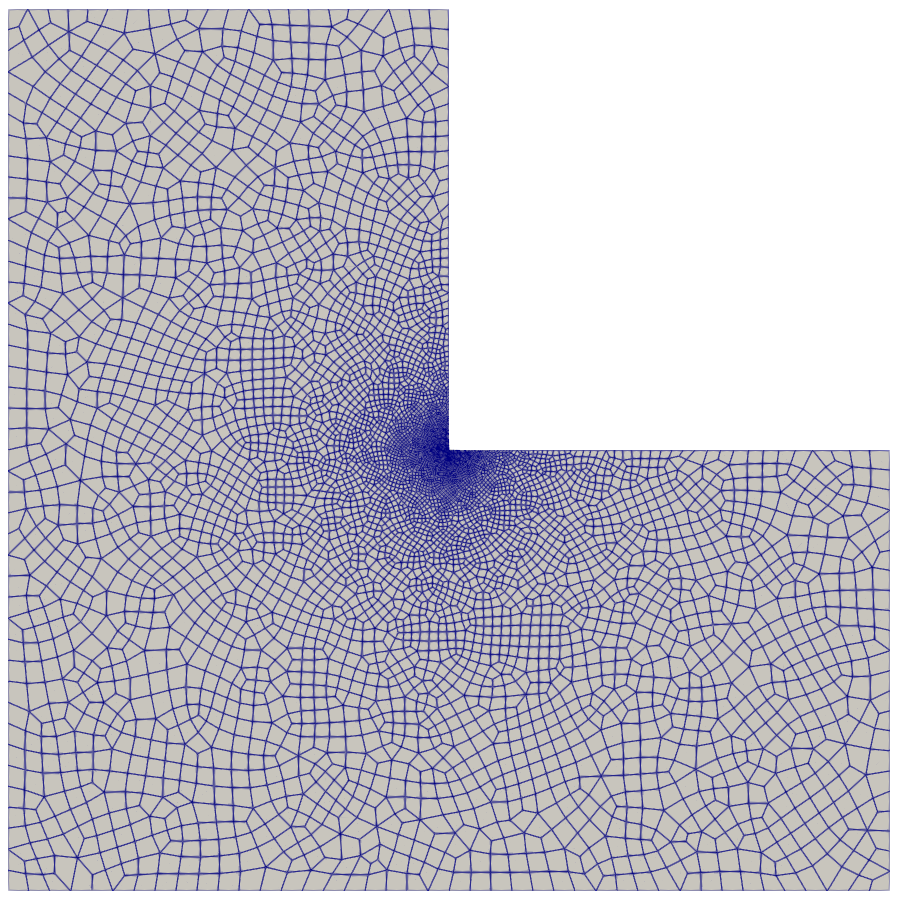}}\qquad%
    \caption{Four of the levels employed for the L-shaped domain, refined near the re-entrant corner.}
    \label{fig:levels_L_shape}
\end{figure}

A sequence of unstructured L-shaped domains is constructed using \textsc{gmsh}. The right-hand side and 
boundary conditions force a singular behavior at the re-entrant corner. For this reason, several refinements
have been applied near the corner during the generation of the levels. The number of required iterations and the number of degrees of freedom
on the finest level are reported in Table \ref{tab:L_shape}. For every level and polynomial degree, an
almost constant number of iterations required by the solver are observed. A set of the first four levels is shown in Figure \ref{fig:levels_L_shape}.
To highlight the flexibility of the present framework, we repeat the same test with the exception that the coarsest
level is a classical structured L-shaped domain with $192$ cells, keeping all the other levels unchanged. The resulting iteration counts are reported in the
second columns of Table~\ref{tab:L_shape}, indicating robust convergence also in that scenario.

\begingroup
\renewcommand{\arraystretch}{1.1}
\begin{center}
    \begin{table}[!t]
        \centering 
        \begin{tabular}{l|ccr|ccr|ccr|ccr}
            \toprule
            \emph{l}& \multicolumn{3}{|c}{$p=1$} & \multicolumn{3}{|c}{$p=2$} & \multicolumn{3}{|c}{$p=3$} & \multicolumn{3}{|c}{$p=4$}  \\
            \cmidrule(lr){2-4}\cmidrule(lr){5-7}\cmidrule(lr){8-10}\cmidrule(lr){11-13}
            &\emph{\#i} & \emph{\#i\_s} & \#DoFs &\emph{\#i} & \emph{\#i\_s} & \#DoFs &\emph{\#i} & \emph{\#i\_s} & \#DoFs &\emph{\#i} & \emph{\#i\_s} & \#DoFs  \\
            \midrule
            2 & 4 & 5 & 1 691  & 6 & 6 & 6 613  & 9 & 7 & 14 767   & 12 & 11 & 26 153  \\
            3 & 4 & 4 & 2 988  & 6 & 5 & 11 757 & 9 & 7 & 26 308   & 12 & 11 & 46 641  \\
            4 & 4 & 4 & 5 055  & 6 & 5 & 19 965 & 9 & 6 & 44 731   & 12 & 11 & 79 353  \\
            5 & 4 & 4 & 8 735  & 6 & 5 & 34 605 & 9 & 6 & 77 611   & 12 & 11 & 137 753  \\
            6 & 4 & 4 & 14 675 & 6 & 5 & 58 261 & 9 & 5 & 130 759  & 11 & 10 & 232 169  \\
            7 & 4 & 4 & 24 308 & 6 & 5 & 96 661 & 9 & 6 & 217 060  & 12 & 11 & 385 505  \\
            \bottomrule
        \end{tabular}
        \caption{Number of required iterations and DoFs per level to solve the system with the L-shaped domain. Here the column \emph{i\_s} refers to the iteration counts
        when a coarser structured mesh with $192$ elements is used as first level.}
        \label{tab:L_shape}
    \end{table}
\end{center}
\endgroup

\subsubsection{Fichera's corner}
\begin{figure}[!h]
    \subfloat{\includegraphics[width=0.33\linewidth]{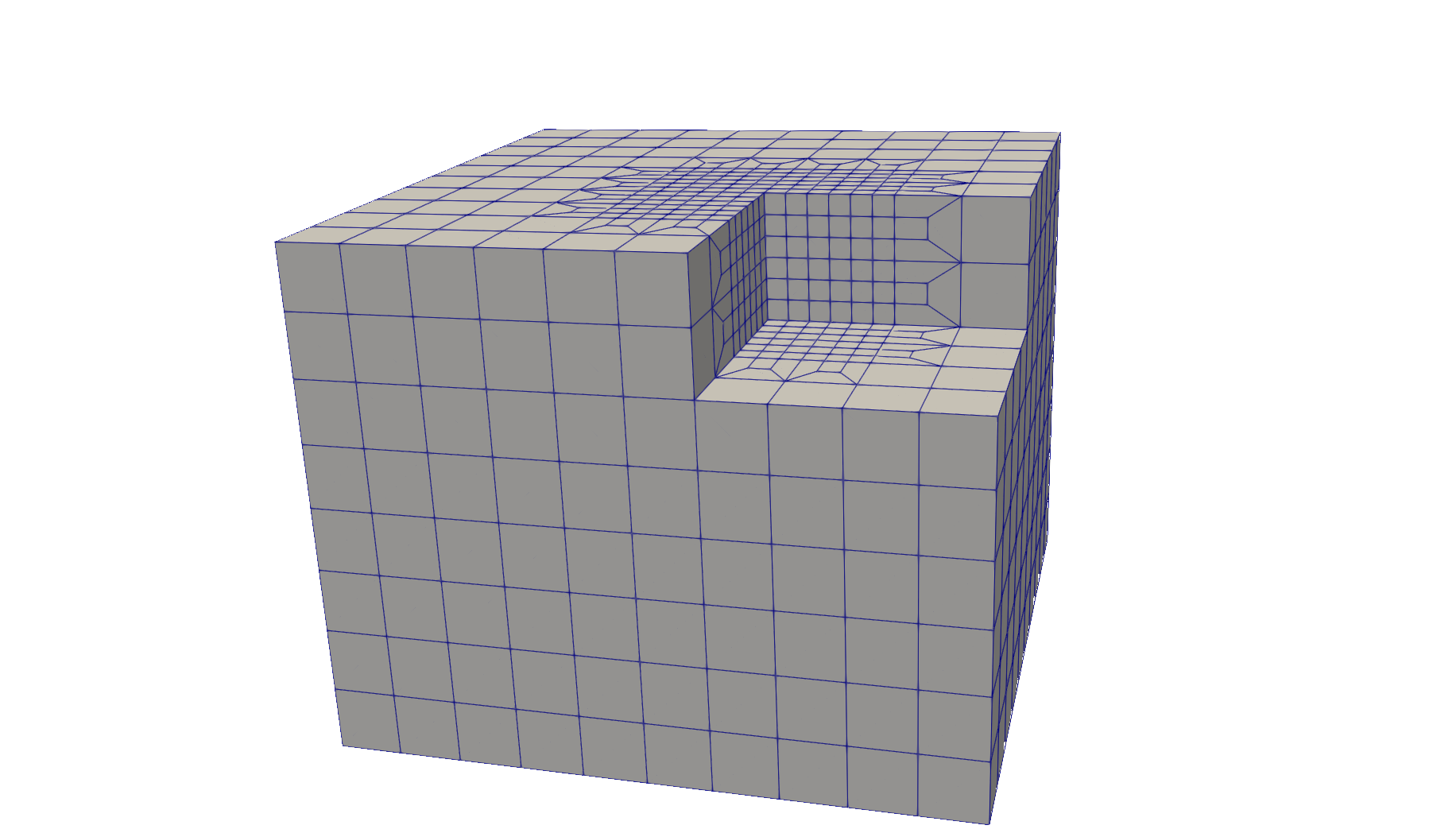}}\qquad%
    \subfloat{\includegraphics[width=0.33\textwidth]{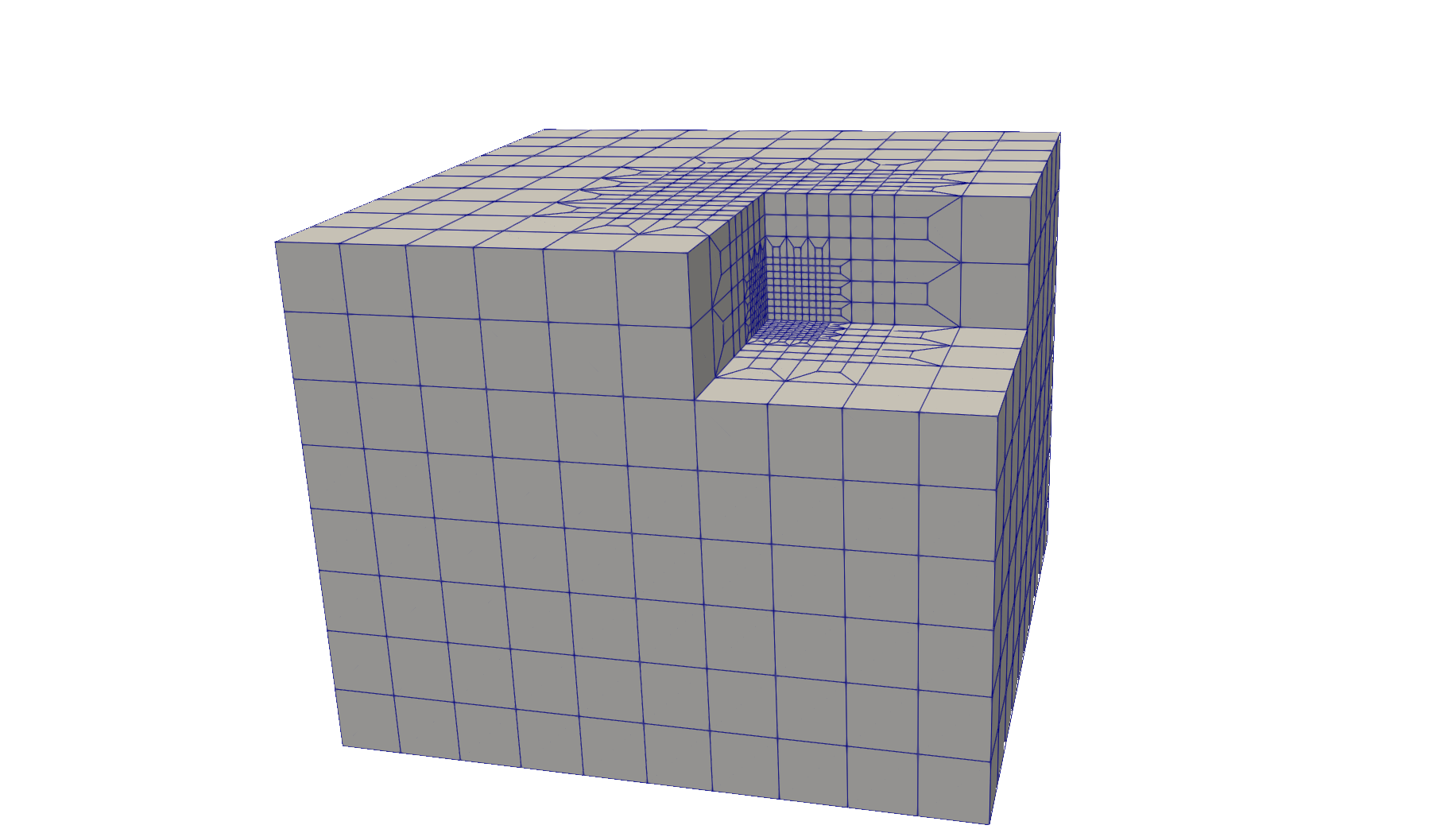}}\qquad%
    \subfloat{\includegraphics[width=0.33\textwidth]{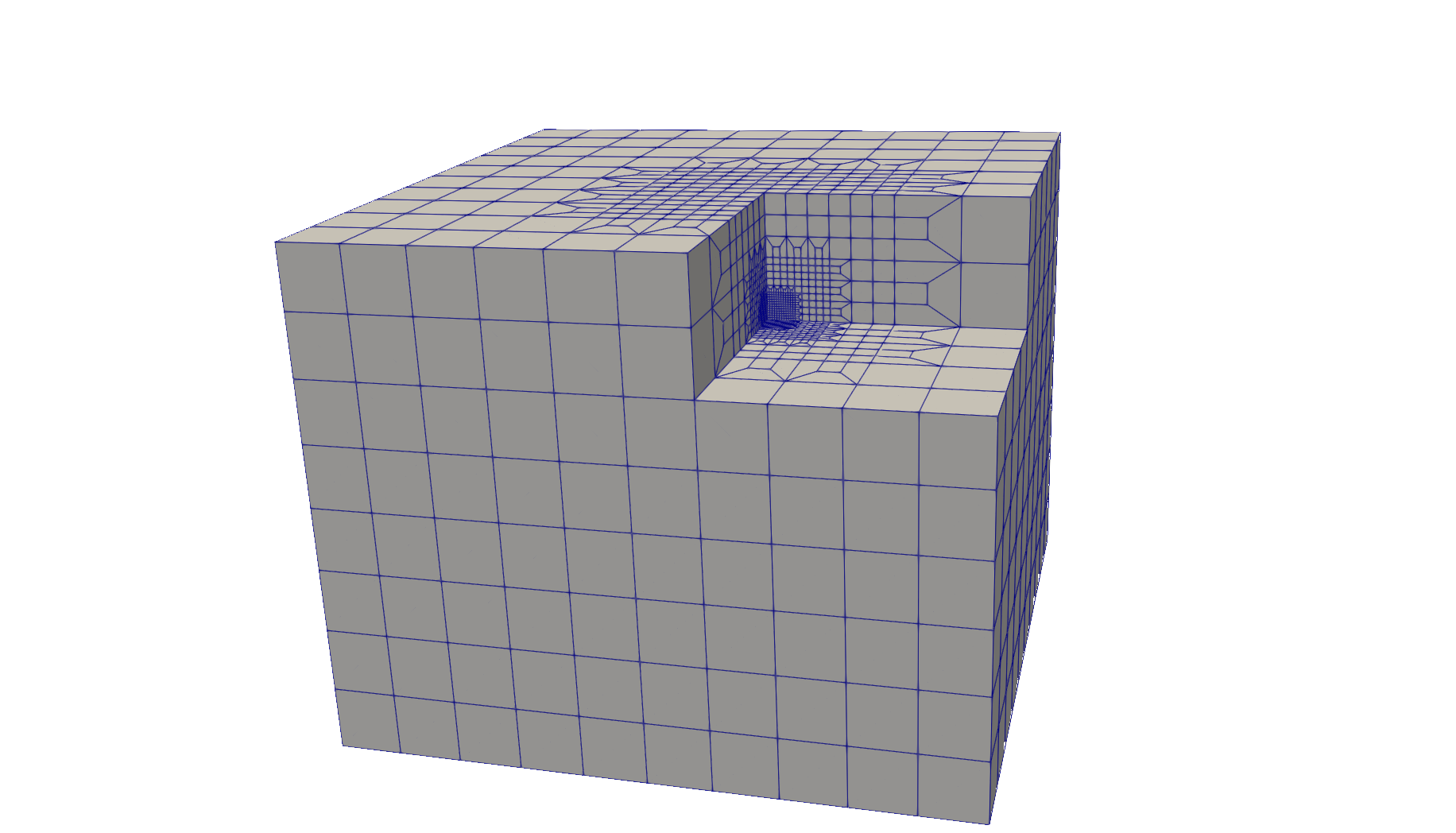}}\qquad%
    \caption{Consecutive levels employed for the Fichera test, refined near the re-entrant corner.}
    \label{fig:levels_Fichera}
\end{figure}

We consider the \emph{Fichera's corner} test case, which represents the three-dimensional extension of
the two-dimensional L-shaped test shown in \ref{subsubsec:sanity_check_2D_l_shaped}. Levels are generated
with \textsc{Coreform Cubit} starting from a coarse, structured three-dimensional L-shaped domain, refining
around the vertex located at the re-entrant corner by increasing the $\mathtt{element\_depth}$ parameter (how
many elements away from the specified vertex are refined) from $2$ to $5$. A sequence of three consecutive levels is
shown in Figure ~\ref{fig:levels_Fichera}. We observe in~Table~\ref{tab:fichera} that, similarly to the two-dimensional case, for each
polynomial degree $p$ we obtain iteration counts which are robust with respect to the number of levels employed in the hierarchy.

\begingroup
\renewcommand{\arraystretch}{1.1}
\begin{center}
    \begin{table}[!h]
        \centering 
        \begin{tabular}{l|cr|cr|cr|cr}
            \toprule
            \emph{l}& \multicolumn{2}{|c}{$p=1$} & \multicolumn{2}{|c}{$p=2$} & \multicolumn{2}{|c}{$p=3$} & \multicolumn{2}{|c}{$p=4$}  \\
            \cmidrule(lr){2-3}\cmidrule(lr){4-5}\cmidrule(lr){6-7}\cmidrule(lr){8-9}
            &\emph{\#i} & \#DoFs &\emph{\#i} & \#DoFs &\emph{\#i} & \#DoFs &\emph{\#i} & \#DoFs  \\
            \midrule
            2 & 4 & 3  983  & 5 & 30 105    & 7 & 99  775    & 9 &  234 401 \\
            3 & 4 & 11 579  & 4 & 90 171    & 7 & 301 789    & 9 &  712 445 \\
            4 & 3 & 27 859  & 4 & 219 283   & 6 & 736 381    & 8 &  1 741 261 \\
            5 & 3 & 57 719  & 4 & 456 513   & 6 & 1 535 311  & 8 &  3 633 041  \\
            6 & 4 & 219 283 & 4 & 1 741 261 & 7 & 5 862 799  & 9 &  13 880 761  \\
            \bottomrule
        \end{tabular}
        \caption{Number of required iterations and DoFs per level to solve the system with the Fichera's corner test.}
        \label{tab:fichera}
    \end{table}
\end{center}
\endgroup

%% file: chapters/applications_fea.tex
\subsection{Applications to FEA}
On a given linear elastic body $\Omega \subset \mathbb{R}^3$, with boundary $\partial \Omega$ partitioned in $\partial \Omega_D$ and $\partial \Omega_N$, we solve (under the assumption of infinitesimal strains) the linear elasticity equation for compressible materials. The governing equations for the unknown displacement $\boldsymbol{u}: \Omega \rightarrow \mathbb{R}^3$ are (\cite{Gurtin}):
\begin{equation}
    \label{eqn:linear_elasticity}
        \begin{cases}
        -\nabla\cdot\sigma(\boldsymbol{u}) &= \boldsymbol{f} \qquad \text{in}\ \Omega, \\
        \sigma(\boldsymbol{u}) &= \lambda\,\hbox{tr}\,(\varepsilon(\boldsymbol{u})) I + 2\mu\varepsilon(\boldsymbol{u}),\\ 
        \varepsilon(\boldsymbol{u}) &= \frac{1}{2}\left(\nabla \boldsymbol{u} + (\nabla \boldsymbol{u})^{\top}\right),\\
        \boldsymbol{u} &= \boldsymbol{0} \qquad \text{on}\ \partial \Omega_D,\\
        \sigma(\boldsymbol{u}) \cdot \boldsymbol{n} &= \boldsymbol{g} \qquad \text{on}\ \partial \Omega_N,
    \end{cases}
    \end{equation}
where $\lambda, \mu$ are the Lamé coefficients for the material, $\sigma$ is the stress tensor, $I$ the identity tensor, $\varepsilon(\boldsymbol{u})$ the linearized strain rate tensor and $\boldsymbol{f}: \Omega \rightarrow \mathbb{R}^3$ is the body force exerted per unit volume.
The boundary $\partial \Omega_D$ is clamped, whereas the normal load $\boldsymbol{g}$ is imposed on $\partial \Omega_N$. A natural space for the kinematically admissible displacement field $u$ and related test functions is $V_D \coloneqq \{ \boldsymbol{v} \in [H^1(\Omega)]^3 : \boldsymbol{v}|_{\partial \Omega_D}=\boldsymbol{0} \}$, the set
of vector-valued $H^1$ functions in $\Omega$ with zero trace (displacement) on the Dirichlet portion of the boundary $\partial \Omega_D$. The variational formulation of \eqref{eqn:linear_elasticity} requires to find $\boldsymbol{u} \in V_D$ such that
$$
a(\boldsymbol{u},\boldsymbol{v})=b(\boldsymbol{v}) \qquad \forall \boldsymbol{v} \in V_D,
$$
where $a(\boldsymbol{u},\boldsymbol{v})=\bigl(\sigma(\boldsymbol{u}), \varepsilon({\boldsymbol{v}}) \bigr)_{\Omega}$ and $b(\boldsymbol{v})=\bigl( f, \boldsymbol{v}\bigr)_{\Omega} + \langle \boldsymbol{g},\boldsymbol{v} \rangle_{\partial \Omega_N}$. Well-posedness for this problem follows from Korn inequalities
and the Lax-Milgram lemma \cite{ErnGuermond}. The material considered in the forthcoming examples is plain steel with Young's modulus $E=205$ \si{GPa} and Poisson's ratio $\nu=0.3$. Given $E$ and $\nu$, Lamè constants can be computed as follows:
\begin{align}
    \lambda &= \frac{E \nu}{(1-2\nu)(1+\nu)},\\
    \mu &= \frac{E}{2(1+\nu)}.
\end{align}
It is well known that low-order elements ($p=1,2$) applied to \eqref{eqn:linear_elasticity} with mixed boundary conditions may suffer from the so-called \emph{locking} phenomenon when $\lambda \rightarrow \frac{1}{2}$, meaning that the material
is nearly incompressible. Among the possible remedies, a simple and generally accepted rule in engineering is that using higher order conforming elements can reduce the potential for locking (even though there are counterexamples,
see \cite{AINSWORTH2022115034}). On the other hand, matrix-vector multiplications for higher order matrix-based discretization rapidly become quite expensive as the bandwidth of the matrices increases \cite{KronbichlerPerssonDG, Davydov2020}. Other
strategies that can prevent locking consist in mixed or augmented formulations \cite{bbf}, but higher order polynomials are rarely employed in practice, usually owing to some low regularity arguments. The forthcoming examples consider
domains arising from applications in structural analysis that are widely studied in literature. The geometry is described through $\mathtt{.step}$ files, a standardized ISO file format used in CAD design. Models have been first repaired and
later meshed using \textsc{Coreform Cubit}. On both domains, we solve the elasticity equations varying the polynomial degree $p$, the number of levels in the hierarchy and the number of processes. We show in Table \ref{tab:tests_fea} the chosen solvers depending on the degree of the finite element space.
For $\mathcal{Q}^1$ elements we compare with AMG, a setting very competitive for linear elements. For higher orders, we compare with polynomial multigrid (PMG) preconditioning, using CG preconditioned by AMG as a coarse grid solver since the coarsest mesh of the polynomial hierarchy is made by linear elements.

\begin{table}[!h]
    \begin{center}
        \centering
        \begin{tabular}{|l|r|r|r|r|}%
     \hline
     \multicolumn{4}{|c|}{\textbf{Solvers for problem \eqref{eqn:linear_elasticity}} } \\
     \hline
     \diagbox{Element}{Solver} & AMG & NN & PMG  \\ 
     \hline
      $\mathcal{Q}^1$ &\cmark  &\cmark &\xmark \\
      $\mathcal{Q}^2$ &\cmark  &\cmark &\cmark  \\
      $\mathcal{Q}^3$ &\xmark  &\cmark &\cmark  \\
      $\mathcal{Q}^4$ &\xmark  &\cmark &\cmark  \\
     \hline
    \end{tabular}
\end{center}
\caption{In each row we report the polynomial space and the solvers applied to the FEA problem \eqref{eqn:linear_elasticity}. Legend: AMG (Algebraic multigrid), NN (Non-nested multigrid), PMG (Polynomial multigrid).}
\label{tab:tests_fea}
\end{table}
In order to successfully apply AMG to vector-valued problems it is necessary to provide the nullspace of the weak form $\tilde{a}(\cdot,\cdot)$ associated to the problem with free boundary conditions on $\partial \Omega$, defined as:
\begin{equation}
    \label{eqn:kernel_a}
    \ker(\tilde{a}) \coloneqq \{ \bm{v} \in H^1(\Omega): \tilde{a}(\bm{v},\bm{v})=0 \}.
\end{equation}
In case of the three-dimensional linear elasticity problem with strain tensor $\varepsilon$ it holds~\cite{AMG_linear_elasticity}:
\begin{equation}
    \label{eqn:rbm}
    \ker(\tilde{a}) = \ker(\varepsilon)= \spn\{ \bm{t}_1,\bm{t}_2,\bm{t}_3,\bm{r}_4,\bm{r}_5,\bm{r}_6 \},
\end{equation}
with
\begin{align}
    \bm{t}_1 = \bm{e}_1, \bm{t}_2 = \bm{e}_2,  \bm{t}_3 = \bm{e}_3, \bm{r}_4 = \begin{bmatrix} 0 \\ z \\ -y \end{bmatrix} , \bm{r}_5 = \begin{bmatrix} -z \\ 0 \\ x \end{bmatrix} , \bm{r}_6 = \begin{bmatrix} y \\ -x \\ 0 \end{bmatrix}.
\end{align}
The first three elements correspond to translation modes and can be represented by the canonical basis $\{\bm{e}_i\}_{i=1}^3$ of $\mathbb{R}^3$, while $\{ \bm{r}_i (\bm{x})\}_{i=4}^6$ are the rotational modes describing rotations around coordinate axes. It should be noted that, from a user's perspective, the human time required for preparing and exporting all levels might certainly be non-negligible in practice. On the other hand, classical AMG implementations such as TrilinosML~\cite{TrilinosML}, the algorithm one used in this work, usually require only the setup of a list of parameters related to the problem at hand.

\subsubsection{Piston}
\begin{figure}[H]
    \centering
    \subfloat[a][CAD model for the piston.]{\includegraphics[width=0.44\textwidth,valign=c]{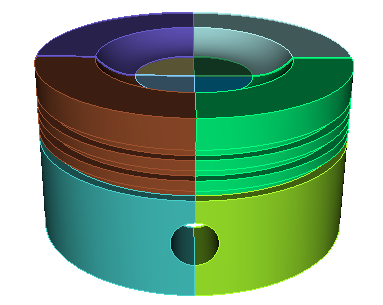}}
    \hfill
    \subfloat[b][Two non-nested levels for the piston (clipped, reference configuration).]{\includegraphics[width=0.49\textwidth,valign=c]{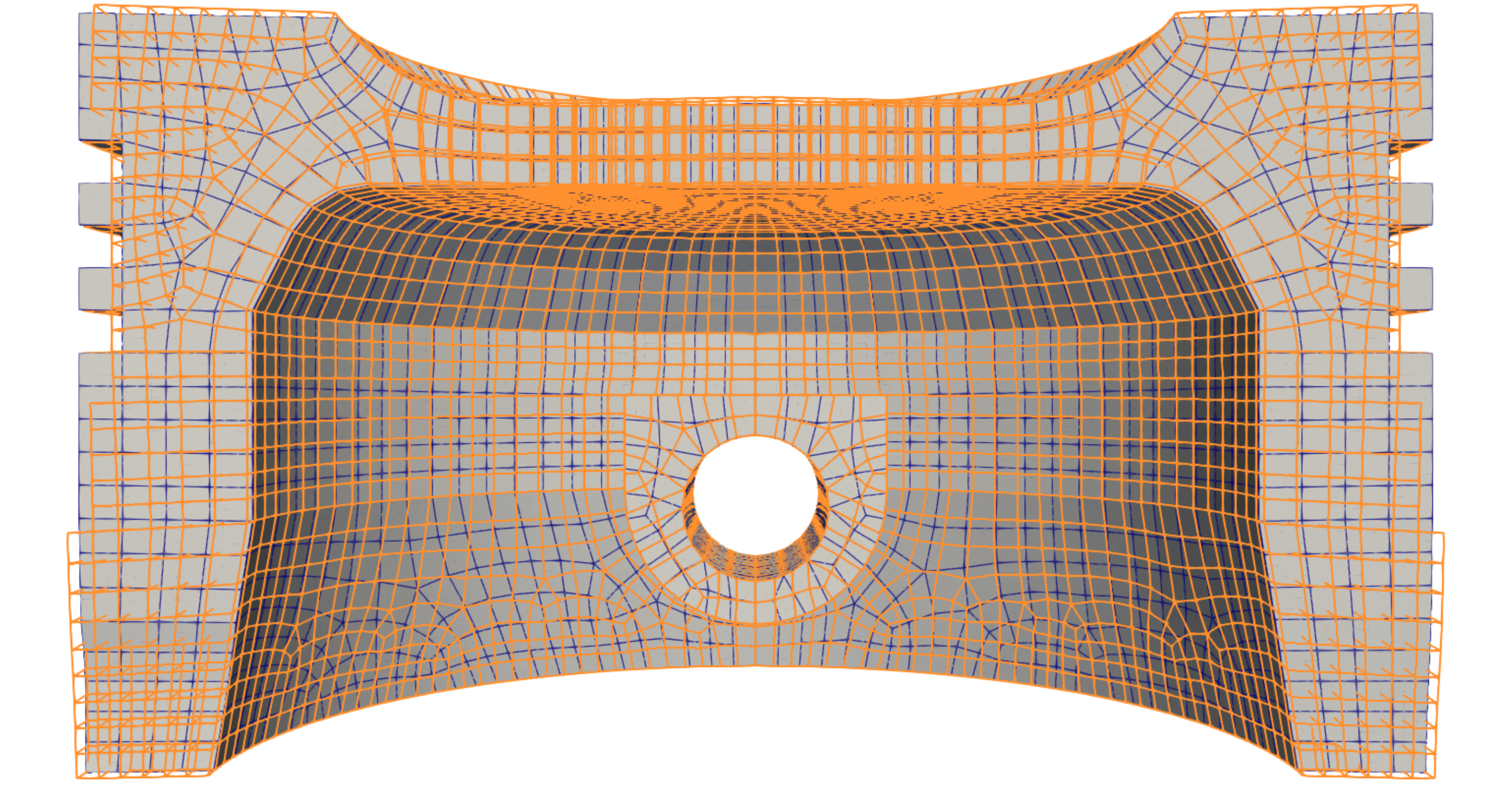}}
    \caption{Left: CAD model used in the preprocessing procedure. Right: Two different levels. For the sake of visualization, the finer level is displayed in orange using a wireframe representation.
        The coarser level is represented as a volumetric mesh (edges in blue). Notice that each level has not been generated on top of the coarser one through a global or local refinement process.}
        \label{fig:piston}
  \end{figure}

\begin{figure}[H]
    \centering
    \includegraphics[width=6cm]{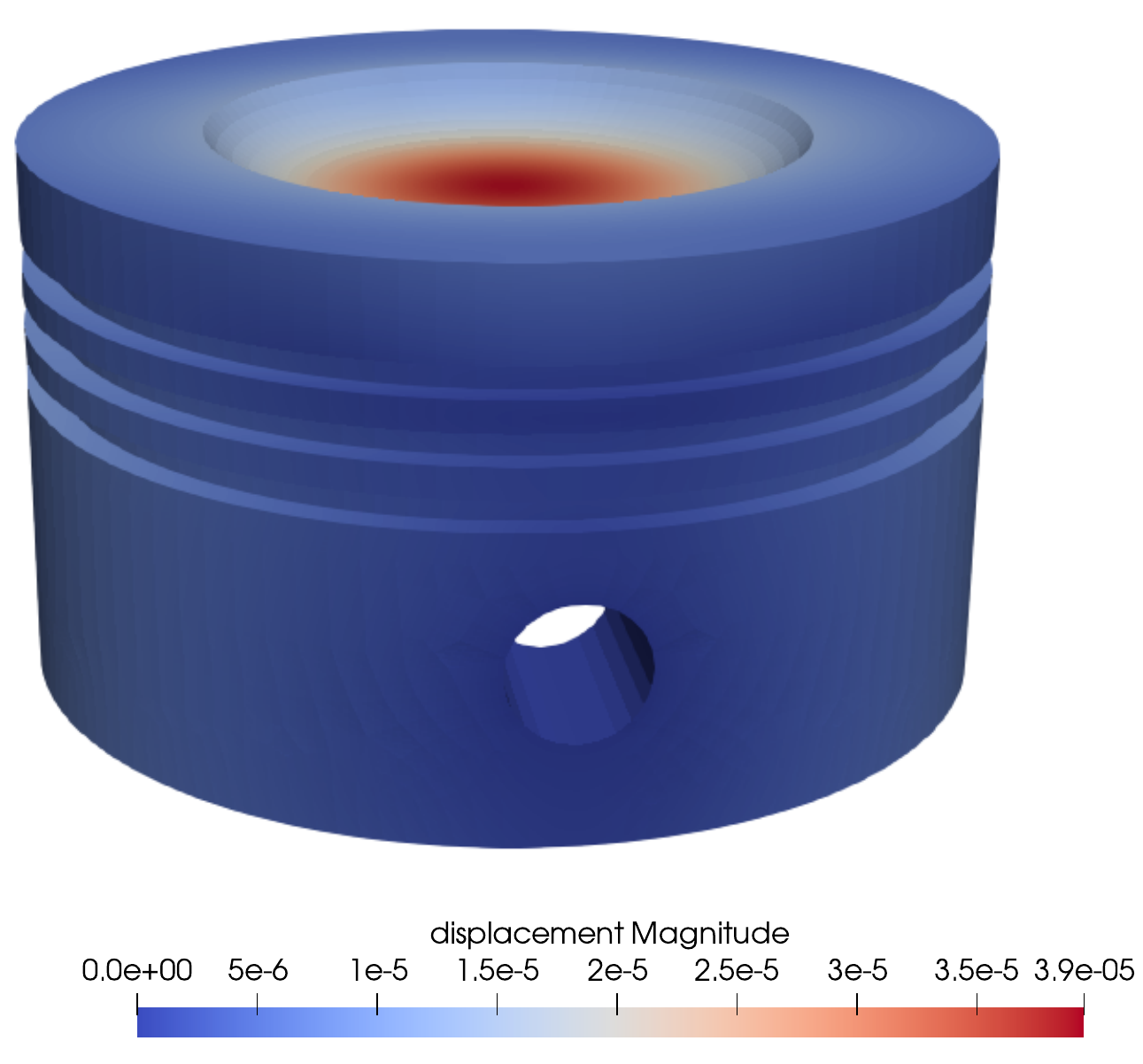}
    \hfill
    \includegraphics[width=6cm]{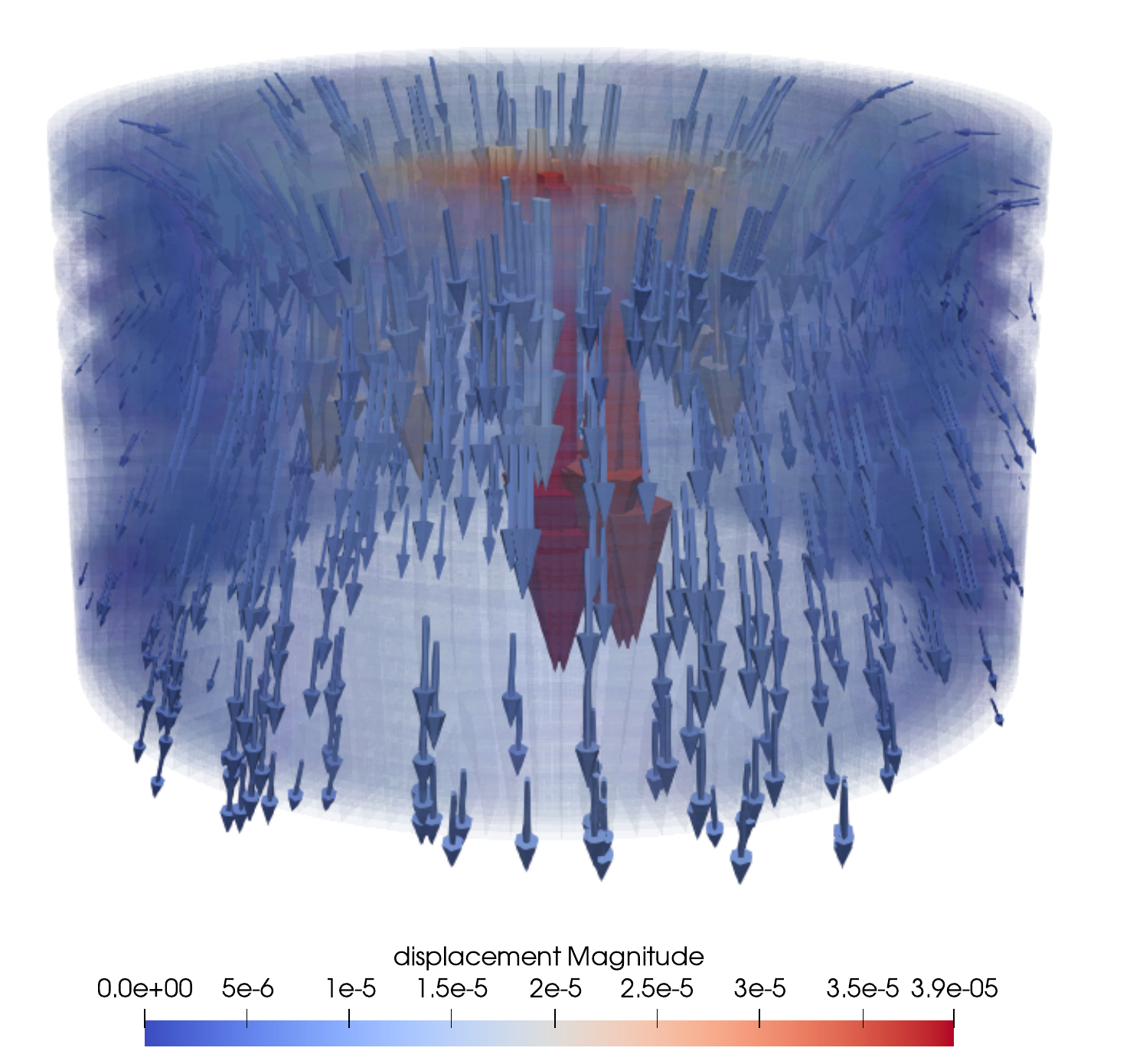}
    \caption{Left: Magnitude of the displacement vector $\boldsymbol{u}$ for the piston test case. Right: Scaled view of the vector field $u$.}
    \label{fig:piston_sols}
\end{figure}

The domain depicted in Figure \ref{fig:piston}(a) describes a flat head piston. We set $\boldsymbol{f} = \boldsymbol{0}$ $\si{N/m^3}$, neglecting external body forces, and we prescribe zero displacement on the surfaces of pin hole and a normal load $\boldsymbol{g} = -10^5 \boldsymbol{n}$ $\si{ Pa}$ on the head, being $\boldsymbol{n}$ the outward unit normal. The rest of the boundary is traction-free.
Table~\ref{tab:DoFs_per_meshsize_piston} lists the approximate mesh size per level $h_l$ set inside the mesh generator and the resulting number of DoF obtained using continuous Lagrangian finite elements for each $l \in \{1,\ldots,4\}$ while varying the polynomial degree $p$ from 1 to 4. Each generated grid constitutes a level used in the V-cycle. Figure \ref{fig:piston}(b) shows a zoom to highlight the non-nestedness of two resulting levels obtained by the aforementioned process. A visualization of the magnitude of the displacement $\bm u$ is given in Figure \ref{fig:piston_sols}. Table~\ref{tab:piston} shows time to solution and number of iterations for polynomial degrees $p$ ranging from $1$ to $4$ and different number of levels, with
$48$ processes. Some general observations can be made. First, the number of iterations for the non-nested multigrid method is roughly constant for every degree $p$, indicating a correct implementation of the proposed methodology. Moreover, for a fixed number of levels $l$, we observe similar patterns in the time to solution. Figure~\ref{fig:times_piston} shows the times to solution for different solvers across the number of levels when $p=2$. For $\mathcal{Q}^1$ elements, AMG shows higher number of iterations compared to the non-nested multigrid algorithm, but with overall better time to solution.
When employing $\mathcal{Q}^2$ elements, the non-nested approach shows better times to solution compared to both AMG and PMG. Moving to the $\mathcal{Q}^3$ case, we note comparable results for PMG and the non-nested procedure, while the gap between them starts increasing for higher order elements such as $\mathcal{Q}^4$. This is expected as the prolongation and restriction operators for polynomial global coarsening are employing matrix-free kernels based on sum-factorization
algorithms for classical nodal FEM spaces, which generally have favorable complexity for higher-order elements and tensor-product quadrature rules.

\begingroup
\renewcommand{\arraystretch}{1.1}
    \begin{table}[!t]
        \centering
        \resizebox{\textwidth}{!}{
            \begin{tabular}{l|cc|cc|cc|cc|cc|cc|cc|cc|cc}
            \toprule
            \multicolumn{19}{c}{$48$ processes}  \\
            \toprule
            \multicolumn{5}{c|}{$p=1$} & \multicolumn{6}{c|}{$p=2$} & \multicolumn{4}{c|}{$p=3$} & \multicolumn{4}{c}{$p=4$}  \\
            \toprule
            \emph{l} & \multicolumn{2}{c|}{AMG} & \multicolumn{2}{c|}{NN} & \multicolumn{2}{c|}{AMG} & \multicolumn{2}{c|}{NN} & \multicolumn{2}{c|}{PMG} &  \multicolumn{2}{c|}{NN} & \multicolumn{2}{c|}{PMG} & \multicolumn{2}{c|}{NN} & \multicolumn{2}{c}{PMG} \\
            \cmidrule(lr){2-3}\cmidrule(lr){4-5}\cmidrule(lr){6-7}\cmidrule(lr){8-9}\cmidrule(lr){10-11}\cmidrule(lr){12-13}\cmidrule(lr){14-15}\cmidrule(lr){16-17}\cmidrule(lr){18-19}
            & \emph{\#i} & \emph{t}[s] & \emph{\#i} & \emph{t}[s] &  \emph{\#i} & \emph{t}[s] &  \emph{\#i} & \emph{t}[s] & \emph{\#i} & \emph{t}[s] & \emph{\#i} & \emph{t}[s] & \emph{\#i} & \emph{t}[s] & \emph{\#i} & \emph{t}[s] & \emph{\#i} & \emph{t}[s]   \\
            \midrule
            \input data/piston_data_new
            \bottomrule
        \end{tabular}
        }
    \caption{Number of iterations and time to solution for algebraic multigrid (AMG), polynomial multigrid (PMG), non-nested Multigrid (NN) applied to the piston test case with different polynomial degrees from $p=1$ to $p=4$. AMG times are shown for $\mathcal{Q}^1$ and $\mathcal{Q}^2$ elements only.}
    \label{tab:piston}
    \end{table}
\endgroup

\subsubsection{Wrench}
\begin{figure}[H]
    \centering
    \includegraphics[width=0.8\textwidth]{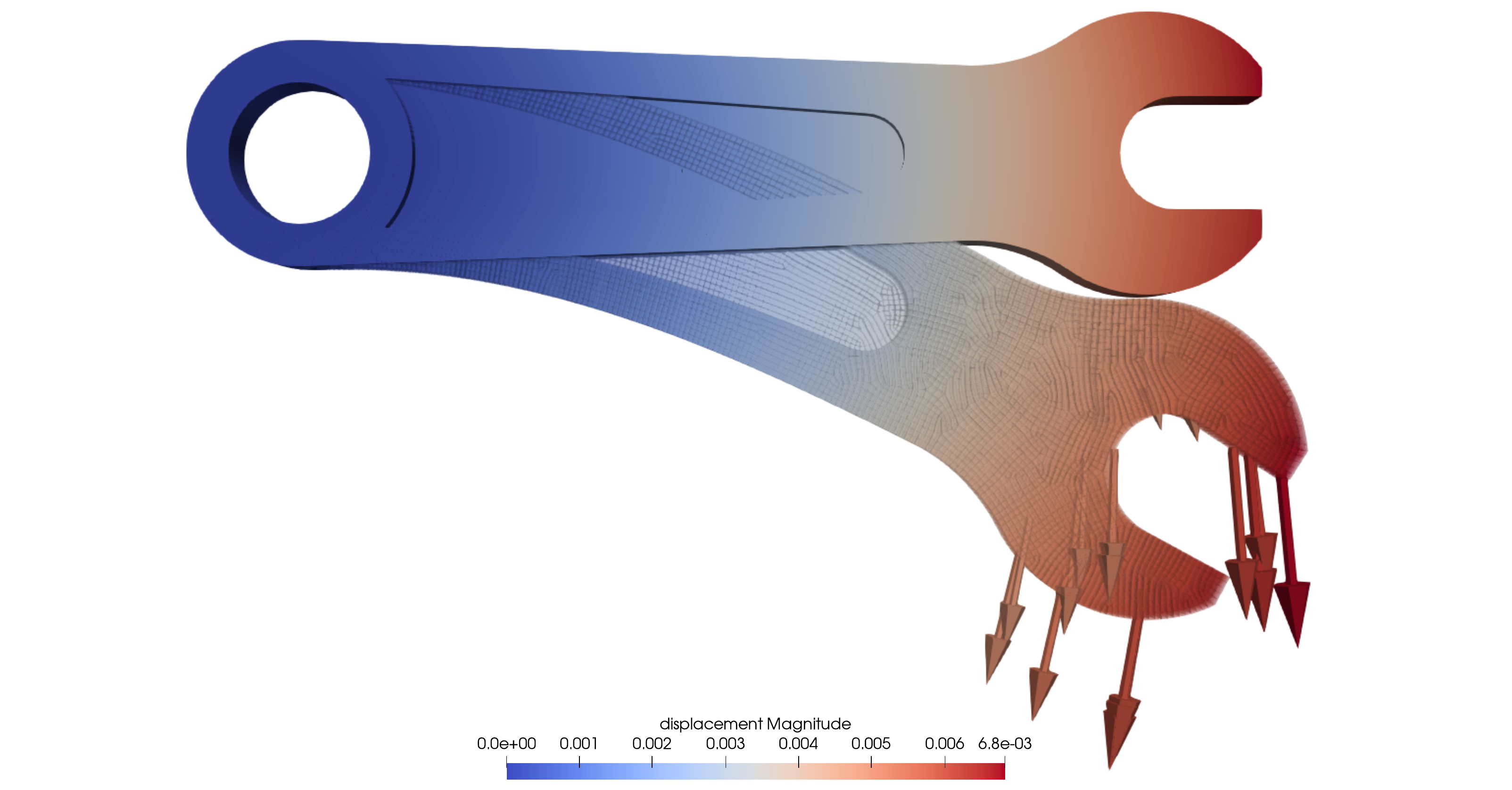}
    \caption{Magnitude (scaled) of the displacement $\boldsymbol{u}$ for the wrench test case.}
    \label{fig:wrench}
\end{figure}

As a second example, we consider the static structural analysis of a wrench, a classical benchmark problem in FEA. Zero displacement is imposed on one head, while a pressure of $\boldsymbol{g} = -10^5 \boldsymbol{n}$ $\si{Pa}$
is applied on one top of the other head, acting in onward direction. The rest of the boundary is traction-free. In Table \ref{tab:wrench} we show time to solution along with the number of required iterations by the outer solver,
while Table \ref{tab:DoFs_per_meshsize_wrench} shows the number of degrees of freedom per each level $l$. The magnitude of the displacement field $\boldsymbol{u}$, along with a graphical representation of the displacement vector field is displayed in Figure \ref{fig:wrench}. The same observations made for
the piston test readily apply also for this test. In particular, Figure~\ref{fig:times_wrench} shows the same behaviour regarding time to solution and polynomial degrees as the one exhibited by the piston test case in Figure~\ref{fig:times_piston}. In this case, we notice a high number of iterations
for AMG with $\mathcal{Q}^1$ and $\mathcal{Q}^2$ elements, while geometric multigrid approaches keep low and constant iteration counts for every polynomial degree $p$, resulting in overall better performances.

\begingroup
\renewcommand{\arraystretch}{1.1}
    \begin{table}[!t]
        \centering
        \resizebox{\textwidth}{!}{
            \begin{tabular}{l|cc|cc|cc|cc|cc|cc|cc|cc|cc}
            \toprule
            \multicolumn{19}{c}{$48$ processes}  \\
            \toprule
            \multicolumn{5}{c|}{$p=1$} & \multicolumn{6}{c|}{$p=2$} & \multicolumn{4}{c|}{$p=3$} & \multicolumn{4}{c}{$p=4$}  \\
            \toprule
            \emph{l} & \multicolumn{2}{c|}{AMG} & \multicolumn{2}{c|}{NN} & \multicolumn{2}{c|}{AMG} & \multicolumn{2}{c|}{NN} & \multicolumn{2}{c|}{PMG} &  \multicolumn{2}{c|}{NN} & \multicolumn{2}{c|}{PMG} & \multicolumn{2}{c|}{NN} & \multicolumn{2}{c}{PMG} \\
            \cmidrule(lr){2-3}\cmidrule(lr){4-5}\cmidrule(lr){6-7}\cmidrule(lr){8-9}\cmidrule(lr){10-11}\cmidrule(lr){12-13}\cmidrule(lr){14-15}\cmidrule(lr){16-17}\cmidrule(lr){18-19}
            & \emph{\#i} & \emph{t}[s] & \emph{\#i} & \emph{t}[s] &  \emph{\#i} & \emph{t}[s] &  \emph{\#i} & \emph{t}[s] & \emph{\#i} & \emph{t}[s] & \emph{\#i} & \emph{t}[s] & \emph{\#i} & \emph{t}[s] & \emph{\#i} & \emph{t}[s] & \emph{\#i} & \emph{t}[s]   \\
            \midrule
            \input data/wrench_data_new
            \bottomrule
        \end{tabular}
        }
    \caption{Number of iterations and time to solution for algebraic multigrid (AMG), polynomial multigrid (PMG), non-nested Multigrid (NN) applied to the wrench test case with different polynomial degrees from $p=1$ to $p=4$. AMG times are shown for $\mathcal{Q}^1$ and $\mathcal{Q}^2$ elements only.}
    \label{tab:wrench}
    \end{table}
\endgroup

\begin{table}[!t]
    \centering
    \begin{minipage}{0.49\textwidth}
        \centering
        \resizebox{\linewidth}{!}{%
        \begin{tabular}{|l|r|r|r|r|r|}
            \hline
            \multicolumn{5}{|c|}{\textbf{DoFs per mesh size and polynomial degree $p$}} \\
            \hline
            \diagbox{$h_l$}{$p$} & 1 & 2 & 3 & 4  \\ %
            \hline
            0.25  &  25 005  &   171 765 &   548 367 & 1 262 901 \\
            0.19  &  50 814  &   360 843 & 1 167 090 & 2 706 561 \\
            0.15  & 110 904  &   804 399 & 2 622 462 & 6 107 073 \\
            0.123 & 165 069  & 1 208 187 & 3 952 203 & 9 219 969 \\
            \hline
        \end{tabular}
        }
        \caption{Total number of degrees of freedom for polynomial degree $p$ and mesh size for the piston test.}
        \label{tab:DoFs_per_meshsize_piston}
    \end{minipage}
    \hfill
    \begin{minipage}{0.49\textwidth}
        \centering
        \resizebox{\linewidth}{!}{%
        \begin{tabular}{|l|r|r|r|r|r|}
            \hline
            \multicolumn{5}{|c|}{\textbf{DoFs per mesh size and polynomial degree $p$}} \\
            \hline
            \diagbox{$h_l$}{$p$} & 1 & 2 & 3 & 4  \\ %
            \hline
            1.8   &  25 425  &   171 138 &   542 691 & 1 245 636 \\
            1.5   &  48 366  &   336 822 & 1 082 052 & 2 500 740 \\
            1.1   & 108 492  &   777 894 & 2 526 372 & 5 872 092 \\
            0.95  & 175 020  & 1 275 312 & 4 166 334 & 9 713 544 \\
            \hline
        \end{tabular}
        }
        \caption{Total number of degrees of freedom for polynomial degree $p$ and mesh size for the wrench test.}
        \label{tab:DoFs_per_meshsize_wrench}
    \end{minipage}
\end{table}

\begin{figure}[h!]
    \centering
    \begin{minipage}{0.49\textwidth}
        \centering
        \begin{tikzpicture}
            \begin{axis}[
                height=7cm,
                width=\textwidth,
                grid=both,
                xlabel={Levels (\( l \))},
                ylabel={Time [s]},
                legend pos=north west,
                xtick={2,3,4},
                scaled y ticks = false,
                legend style={font=\small}
            ]
                \addplot coordinates {(2, 0.5314) (3, 1.1318) (4, 1.4998)};
                \addlegendentry{AMG }
    
                \addplot coordinates {(2, 0.1395) (3, 0.2119) (4, 0.3607)};
                \addlegendentry{NN }
    
                \addplot coordinates {(2, 0.2081) (3, 0.3835) (4, 0.5841)};
                \addlegendentry{PMG }
            \end{axis}
        \end{tikzpicture}
        \caption{Piston test case (\( p=2 \)).}
        \label{fig:times_piston}
    \end{minipage}
    \hfill
    \begin{minipage}{0.49\textwidth}
        \centering
        \begin{tikzpicture}
            \begin{axis}[
                height=7cm,
                width=\textwidth,
                grid=both,
                xlabel={Levels (\( l \))},
                legend pos=north west,
                xtick={2,3,4},
                scaled y ticks = false,
                legend style={font=\small}
            ]
                \addplot coordinates {(2, 0.5947) (3, 1.5574) (4, 3.2584)};
                \addlegendentry{AMG }
                
                \addplot coordinates {(2, 0.2381) (3, 0.3601) (4, 0.5777)};
                \addlegendentry{NN }
                
                \addplot coordinates {(2, 0.1889) (3, 0.5331) (4, 0.8170)};
                \addlegendentry{PMG }
            \end{axis}
        \end{tikzpicture}
        \caption{Wrench test case (\( p=2 \)).}
        \label{fig:times_wrench}
    \end{minipage}
    \caption{Time to solution for different solvers across levels (\( l \)) for \( p=2 \) in piston and wrench test cases.}
    \label{fig:comparison}
\end{figure}

%% file: data/piston_data_new.tex
2 & 17 & 0.0279 & 9 & 0.0817 & 21 & 0.5314 & 10 & 0.1395 & 8 & 0.2081 & 11 & 0.2519 & 12 & 0.3516 & 12 & 0.5200 & 9  & 0.3641 \\
3 & 16 & 0.0503 & 9 & 0.1074 & 21 & 1.1318 & 10 & 0.2119 & 8 & 0.3835 & 11 & 0.5057 & 11 & 0.6260 & 13 & 1.3440 & 10 & 0.8391 \\
4 & 15 & 0.0765 & 9 & 0.1369 & 20 & 1.4998 & 11 & 0.3607 & 8 & 0.5841 & 12 & 0.9646 & 12 & 1.0432 & 13 & 2.6388 & 10 & 1.2863 \\

%% file: data/wrench_data_new.tex
2 & 29 & 0.0347 & 13 & 0.1540 & 30 & 0.5947 & 14 & 0.2381 & 6 & 0.1889 & 15 & 0.4192 & 8 & 0.2855 & 16 & 0.8569 & 9 & 0.3469 \\
3 & 33 & 0.0347 & 13 & 0.1932 & 33 & 1.5574 & 14 & 0.3601 & 6 & 0.5331 & 16 & 0.9013 & 8 & 0.7643 & 17 & 2.1006 & 8 & 0.8135 \\
4 & 37 & 0.2035 & 13 & 0.2482 & 40 & 3.2584 & 14 & 0.5777 & 6 & 0.8170 & 16 & 1.6505 & 8 & 1.1902 & 17 & 4.4582 & 8 & 1.6841  \\

%% file: chapters/performance_evaluation.tex
\section{Performance evaluation}
\label{sec:performance}
In the following, the performance of our implementation is investigated by means of the 3D Poisson problem detailed in Section \ref{sec:numerical_experiments}.
In the next tests, we compare the time to solution and the breakdown of times spent on each multigrid level by global coarsening and non-nested multigrid for different geometries for which both methods
can be applied. In order to perform a comparison, we employ identical hierarchies of nested meshes. As a crucial remark, we stress that we are measuring against the highly optimized baseline of matrix-free methods described in \cite{MHPSK}. We consider two types of 3D meshes:
\begin{itemize}
    \item \textbf{cube}: refine globally $l$ times the cube $[-1,1]^3$,
    \item \textbf{ball}: refine globally $l$ times $\mathcal{B}_1(\boldsymbol{0})$, the ball of radius $1$ centered at the origin.
\end{itemize}

\begin{table}[!t]
    \centering
    \begin{minipage}{0.48\textwidth}
        \centering
        \begingroup
        \renewcommand{\arraystretch}{1.1}
        \resizebox{\linewidth}{!}{
        \begin{tabular}{l|cccc|cccc}
            \toprule
            \multicolumn{5}{c}{1 process} & \multicolumn{4}{|c}{12 processes}  \\
            \toprule
            \emph{l} & \multicolumn{2}{c}{GC} & \multicolumn{2}{c}{NN} & \multicolumn{2}{|c}{GC} & \multicolumn{2}{c}{NN}  \\
            \cmidrule(lr){2-3}\cmidrule(lr){4-5}\cmidrule(lr){6-7}\cmidrule(lr){8-9}
            & \emph{\#i} & \emph{t}[s] & \emph{\#i}  & \emph{t}[s] & \emph{\#i} & \emph{t}[s] & \emph{\#i} & \emph{t}[s]  \\
            \midrule
            4 & 4 & \num{2.2e-2} & 4 & \num{3.3e-2} & 4 & \num{4.7e-3} & 4 & \num{7.9e-3} \\
            5 & 4 & \num{1.6e-1} & 4 & \num{2.5e-1} & 4 & \num{1.8e-2} & 4 & \num{3.0e-2} \\
            6 & 4 & \num{1.3e0}  & 4 & \num{2.0e+0} & 4 & \num{1.5e-1} & 4 & \num{2.1e-1} \\
            7 & 4 & \num{1.1e1}  & 4 & \num{1.6e1}  & 4 & \num{1.2e0} & 4 & \num{1.6e0} \\
            \bottomrule
        \end{tabular}
        }
        \caption{Number of iterations and time to solution for Global Coarsening (GC) and Non-Nested Multigrid (NN) applied to the $\mathtt{cube}$ test case with polynomial degree $p=4$.}
        \label{tab:cube}
        \endgroup
    \end{minipage}
    \hfill
    \begin{minipage}{0.48\textwidth}
        \centering
        \begingroup
        \renewcommand{\arraystretch}{1.1}
        \resizebox{\linewidth}{!}{
        \begin{tabular}{l|cccc|cccc}
            \toprule
            \multicolumn{5}{c}{1 process} & \multicolumn{4}{|c}{12 processes}  \\
            \toprule
            \emph{l} & \multicolumn{2}{c}{GC} & \multicolumn{2}{c}{NN} & \multicolumn{2}{|c}{GC} & \multicolumn{2}{c}{NN}  \\
            \cmidrule(lr){2-3}\cmidrule(lr){4-5}\cmidrule(lr){6-7}\cmidrule(lr){8-9}
            & \emph{\#i} & \emph{t}[s] & \emph{\#i}  & \emph{t}[s] & \emph{\#i} & \emph{t}[s] & \emph{\#i} & \emph{t}[s]  \\
            \midrule
            4 & 5 & \num{2.4e-1}  & 5 & \num{4.2e-1}  & 5 & \num{3.3e-2}  & 5 & \num{5.4e-2} \\
            5 & 6 & \num{2.5e+0}  & 6 & \num{3.5e+0}  & 6 & \num{2.8e-1}  & 6 & \num{3.8e-1}  \\
            6 & 6 & \num{2.1e1}   & 6 & \num{3.0e1}   & 6 & \num{2.3e0}   & 6 &  \num{2.9e+0} \\
            \bottomrule
        \end{tabular}
        }
        \caption{Number of iterations and time to solution for Global Coarsening (GC) and Non-Nested Multigrid (NN) applied to the $\mathtt{ball}$ test case with polynomial degree $p=4$.}
        \label{tab:ball}
        \endgroup
    \end{minipage}
\end{table}

\subsection{Nested meshes: comparison with global coarsening}
\begin{figure}[!t]
    \centering
    \resizebox{\columnwidth}{!}{%
    \begin{tikzpicture}
        \pgfplotsset{width=12cm}
        \begin{axis}
            [
            title= \Large{Cube},
            ymin=0,
            ymax=1.8,
            ybar,
            bar width=0.2,
            legend style={at={(0.4,-0.25)}, 
            anchor=north,legend columns=-1},
            ylabel={Time [s]},
            x tick label style={rotate=45, anchor=east, align=center},
            xtick={0,1,2,3,4,5},
            xticklabels={PreSmoother, Residual, Restrictor, CoarseGridSolver, Prolongator, PostSmoother},
            grid=major,
            grid style={line width=.1pt, draw=gray!50},
            legend pos=north east,
            nodes near coords,
            every node near coord/.append style={rotate=45, anchor=west,font=\footnotesize},
            ]      
            \addplot[blue,fill=blue!30,text=black]coordinates {(0, 7.5879e-01) (1, 3.3142e-01) (2, 8.2615e-01) (3,1.2659e-04) (4, 7.6932e-01) (5,1.0865e+00)};
            \addplot[red,fill=red!30,text=black]  coordinates {(0, 7.6368e-01) (1, 3.3078e-01) (2, 8.3724e-02) (3, 1.2803e-04) (4,8.2194e-02) (5,1.0938e+00)};
            \legend{Non-Nested,Global Coarsening}
        \end{axis} 
    \end{tikzpicture}
    \begin{tikzpicture}
        \pgfplotsset{width=12cm}
        \begin{axis}
            [
            title=\Large{Ball},
            ymin=0,
            ymax=1.8,
            ybar,
            bar width=0.2,
            legend style={at={(0.4,-0.25)}, 
            anchor=north,legend columns=-1},
            x tick label style={rotate=45, anchor=east, align=center},
            xtick={0,1,2,3,4,5},
            xticklabels={PreSmoother, Residual, Restrictor, CoarseGridSolver, Prolongator, PostSmoother},
            grid=major,
            grid style={line width=.1pt, draw=gray!50},
            legend pos=north east,
            nodes near coords,  
            every node near coord/.append style={rotate=45, anchor=west,font=\footnotesize},
            ]
            \addplot[blue,fill=blue!30,text=black] coordinates {(0, 9.9236e-01) (1, 4.5525e-01) (2, 7.2641e-01) (3, 7.7100e-04) (4, 6.7317e-01) (5,1.4400e+00)};
            \addplot[red,fill=red!30,text=black]   coordinates {(0, 9.9452e-01) (1, 4.5411e-01) (2, 7.3621e-02) (3, 6.4481e-04) (4, 7.2175e-02) (5,1.4438e+00)};              
            \legend{Non-Nested, Global Coarsening}
        \end{axis} 
    \end{tikzpicture}
    }
    \caption{Serial profiles of a V-cycle. Left: $\mathtt{Cube}$ example with $l=7$ and $p=4$. Right: $\mathtt{ball}$ example with $l=6$ and $p=4$.}
    \label{fig:breakdown_serial}
\end{figure}
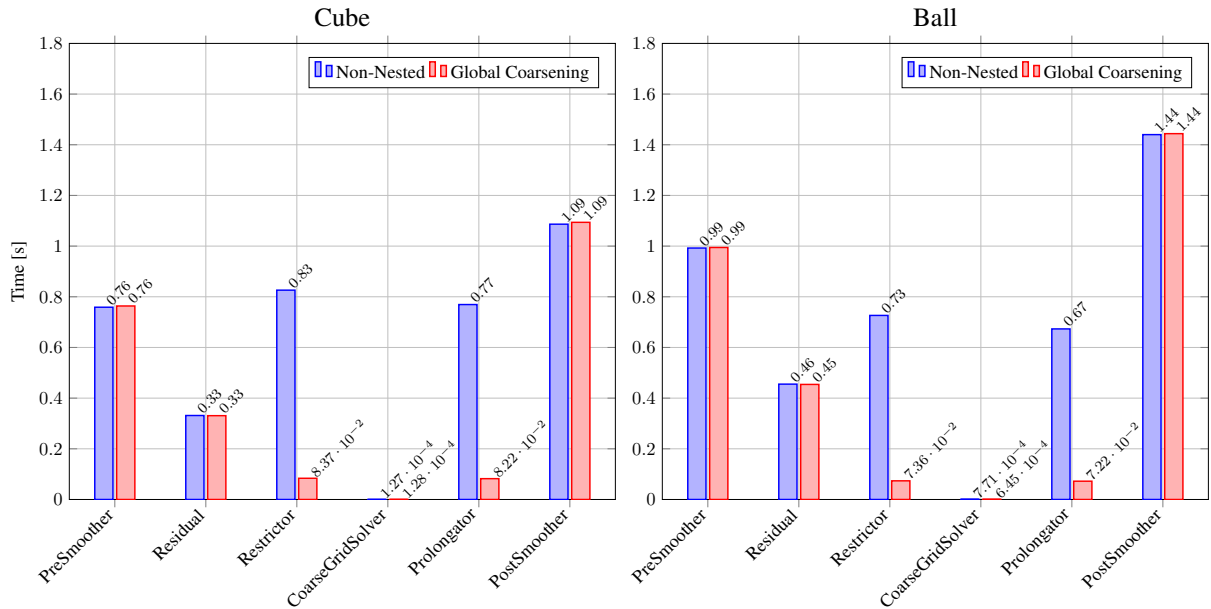

The nestedness of levels for these tests implies that the transfer operator $\Pcf$ for the multigrid method described in \ref{sec:nnmgmethod} coincides with the canonical injection. While the number of iterations is the same for
both variants, the transfers $\Pcf$ and $\mathcal{R}^{(f,c)}\coloneqq \bigl(\Pcf \bigr)^T$ are expected to be much more expensive if compared to classical transfer operators tailored for nested hierarchies.
This is confirmed by Figure~\ref{fig:breakdown_serial} for a serial run with the $\mathtt{cube}$ and $\mathtt{ball}$ examples, where a factor of almost $10$ for intergrid transfers is observed between the two
multigrid variants, whereas the other ingredients are to all extents identical. As shown in Tables \ref{tab:cube} and \ref{tab:ball}, this results in lower times to solution in favor of global coarsening both in the serial
and parallel case.

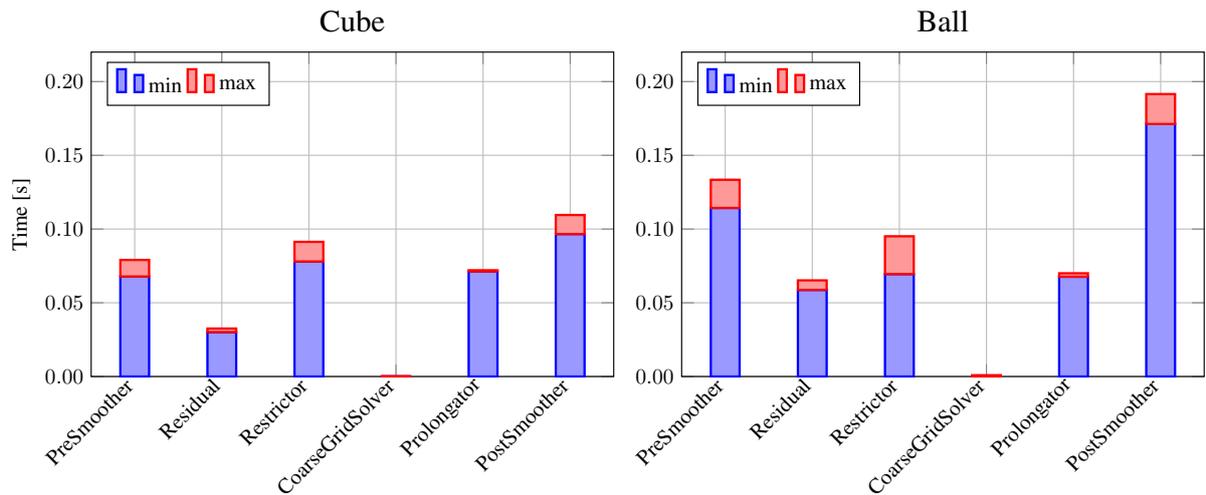
\begin{figure}[!t]
    \centering
    \resizebox{\columnwidth}{!}{%
\begin{tikzpicture}[thick,scale=0.7, every node/.style={scale=0.75}]
    \begin{axis}[
       title=\Large{Cube},
       bar width=10pt,
       width=0.5\textwidth,
       height=0.35\textwidth,
       y tick label style={/pgf/number format/fixed, /pgf/number format/fixed zerofill,/pgf/number format/precision=2},scaled y ticks=false,
       ylabel={Time [s]},
       xmajorgrids,xminorgrids,
       ymajorgrids,yminorgrids,
       ymin=0.0, ymax=0.22,
       xtick = {-1, 0, 1,2,3,4,5,6},
       x tick label style={rotate=45,anchor=east},
       legend columns=3,
       legend pos={north west},
       xticklabels={,PreSmoother, Residual, Restrictor, CoarseGridSolver, Prolongator, PostSmoother}
    ]
    \addlegendimage{ybar,ybar legend,blue, fill=blue!40!white}
    \addlegendentry{min};
    \addlegendimage{ybar,ybar legend,red, fill=red!40!white}
    \addlegendentry{max};
    \addplot [ybar stacked, blue, fill=blue!40!white] coordinates {
    (0, 0.0678343902)
    (1, 0.0301226290)
    (2, 0.0779899665)
    (3, 0.0000284797)
    (4, 0.0712653853)
    (5, 0.0966127082)
    };
    \addplot [ybar stacked, red, fill=red!40!white] coordinates {
    (0,  0.0790534647 - 0.0678343902)
    (1,  0.0325515647 - 0.0301226290)
    (2,  0.0913171172 - 0.0779899665)
    (3,  0.0001322463 - 0.0000284797)
    (4,  0.0720704677 - 0.0712653853)
    (5,  0.1095235770 - 0.0966127082)
    };
    \end{axis}
    \end{tikzpicture} 
    \begin{tikzpicture}[thick,scale=0.7, every node/.style={scale=0.75}]
    \begin{axis}[
       title=\Large{Ball},
       bar width=10pt,
       width=0.5\textwidth,
       height=0.35\textwidth,
       y tick label style={/pgf/number format/fixed, /pgf/number format/fixed zerofill,/pgf/number format/precision=2},scaled y ticks=false,
       xmajorgrids,xminorgrids,
       ymajorgrids,yminorgrids,
       ymin=0.0, ymax=0.22,
       xtick = {-1, 0, 1,2,3,4,5,6},
       x tick label style={rotate=45,anchor=east},
       legend columns=3, legend pos={north west},
       xticklabels={,PreSmoother, Residual, Restrictor, CoarseGridSolver, Prolongator, PostSmoother}
    ]
    \addlegendimage{ybar,ybar legend,blue, fill=blue!40!white}
    \addlegendentry{min};
    \addlegendimage{ybar,ybar legend,red, fill=red!40!white}
    \addlegendentry{max};
    \addplot [ybar stacked, blue, fill=blue!40!white] coordinates {
    (0, 0.1143073748)
    (1, 0.0587237935)
    (2, 0.0694643545)
    (3, 0.0000305518)
    (4, 0.0676649210)
    (5, 0.1712240492)
    };
    \addplot [ybar stacked, red, fill=red!40!white] coordinates {
    (0, 0.1334102550-0.1143073748)
    (1, 0.0651918820-0.0587237935)
    (2, 0.0950921605-0.0694643545)
    (3, 0.0008965638-0.0000305518)
    (4, 0.0700107635-0.0676649210)
    (5, 0.1915411617-0.1712240492)
    };
    \end{axis}
\end{tikzpicture}
}
    \caption{Profile of a V-cycle with 12 processes with non-nested multigrid using a degree 3 Chebyshev smoother. Left: $\mathtt{Cube}$ example with $l=7$ and $p=4$. Right: $\mathtt{ball}$ example with $l=6$ and $p=4$.}
    \label{fig:profile_parallel_NN}
\end{figure}

\noindent Within multigrid literature, a classical assumption is that pre-/post-smoothing steps consume most of the run time, while transfers between levels are, in general, not a bottleneck in the whole pipeline. This holds true whenever nested hierarchies are employed. Indeed, in
Figure~\ref{fig:profile_parallel_NN} we show the minimum and maximum time (measured in seconds) across the parallel processes for each component of the algorithm. The figure illustrates that, for the non-nested approach, transfers are more expensive than the classical transfer
with sum factorization used by global coarsening, with a run time comparable to the application of the smoother.

\begin{figure}[ht]
    \centering
    \resizebox{\columnwidth}{!}{%
    \begin{tikzpicture}[thick,scale=0.7, every node/.style={scale=0.75}]
        \begin{axis}[
           title= \Large{Global Coarsening},
           bar width=10pt,
           width=0.5\textwidth,
           height=0.35\textwidth,
           y tick label style={/pgf/number format/fixed, /pgf/number format/fixed zerofill,/pgf/number format/precision=2},scaled y ticks=false,
           xlabel={Level},
           ylabel={Exclusive time [s]},
           xmajorgrids,xminorgrids,
           ymajorgrids,yminorgrids,
           ymin=0.0, ymax=0.8,
           xtick={0,1,2,3,4,5,6},
           xticklabels={0,1,2,3,4,5,6},
           x tick label style={rotate=45,anchor=east},legend columns=3, legend pos={north west},
        ]
        \addlegendimage{ybar,ybar legend,blue, fill=blue!40!white}
        \addlegendentry{min};
        \addlegendimage{ybar,ybar legend,red, fill=red!40!white}
        \addlegendentry{max};
        
        \addplot [ybar stacked, blue, fill=blue!40!white] coordinates {
            (0, 1.3936e-06)
            (1, 2.29117e-05)
            (2, 2.31637e-05)
            (3, 0.001123314)
            (4, 0.005782064200000001)
            (5, 0.04774718779999999)
            (6, 0.39792133830000)
        };
    
        \addplot [ybar stacked, red, fill=red!40!white] coordinates {
            (0, 0.0002541986 - 1.3936e-06)
            (1,  0.0002491823 - 2.29117e-05)
            (2,  0.0009481009999999999 - 2.31637e-05)
            (3,  0.0023735411000000003 - 0.001123314)
            (4,  0.0072946190000000004 - 0.005782064200000001)
            (5,  0.052839968699999997 - 0.04774718779999999)
            (6,  0.419132155199999 - 0.39792133830000)
        };
        
        \node[above] at (axis cs: 4, 0.007294619) {\num{7.29e-3}};
        \node[above] at (axis cs: 5, 0.052839969) {\num{5.28e-2}};
        \node[above] at (axis cs: 6, 0.419132155) {\num{4.19e-1}};
        \end{axis}
    \end{tikzpicture} 
\begin{tikzpicture}[thick,scale=0.7, every node/.style={scale=0.75}]
    \begin{axis}[
       title= \Large{Non-nested},
       bar width=10pt,
       width=0.5\textwidth,
       height=0.35\textwidth,
       y tick label style={/pgf/number format/fixed, /pgf/number format/fixed zerofill,/pgf/number format/precision=2},scaled y ticks=false,
       xlabel={Level},
       xmajorgrids,xminorgrids,
       ymajorgrids,yminorgrids,
       ymin=0.0, ymax=0.8,
       xtick={0,1,2,3,4,5,6},
       xticklabels={0,1,2,3,4,5,6},
       x tick label style={rotate=45,anchor=east},
       legend columns=3, 
       legend pos={north west},
    ]
    \addlegendimage{ybar,ybar legend,blue, fill=blue!40!white}
    \addlegendentry{min};
    \addlegendimage{ybar,ybar legend,red, fill=red!40!white}
    \addlegendentry{max};
    
    \addplot [ybar stacked, blue, fill=blue!40!white] coordinates {
        (0,  4.96524e-05)
        (1,  0.00035437)
        (2,  0.0008369228)
        (3,  0.0021129927)
        (4,  0.0098132865)
        (5,  0.0745799565)
        (6,  0.60669968)
    };
    
    \addplot [ybar stacked, red, fill=red!40!white] coordinates {
        (0,  0.0002400094 -4.96524e-05     )
        (1,   0.0007251272 - 0.00035437    )
        (2,   0.0012616685 - 0.0008369228    )
        (3,   0.0023541593 - 0.0021129927   )
        (4,   0.010300314100000001 - 0.0098132865   )
        (5,   0.0775995771 - 0.0745799565  )
        (6,   0.625979890 - 0.60669968  )
    };
    
    \node[above] at (axis cs: 4, 0.0103003141) {\num{1.03e-2}};
    \node[above] at (axis cs: 5, 0.0775995771) {\num{7.76e-2}};
    \node[above] at (axis cs: 6, 0.625979890) {\num{6.26e-1}};
    
    \end{axis}
\end{tikzpicture}
}
\caption{Exclusive time in a V-cycle with 12 processes for the $\mathtt{cube}$ example with $l=7$ and $p=4$. Right: global coarsening. Left: non-nested multigrid. Max values have been reported on top of bars associated to
finest levels.}
\label{fig:exclusive_time}
\end{figure}
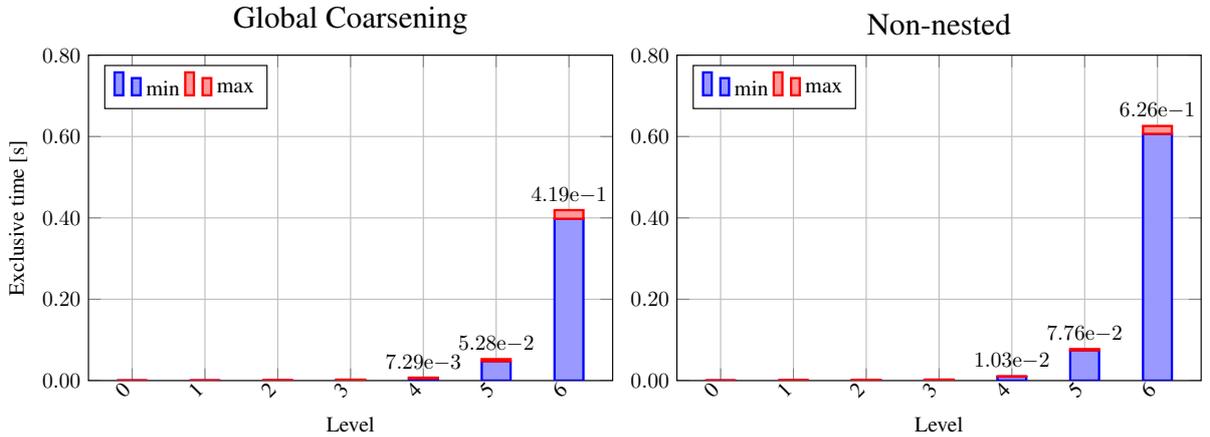

\noindent To gain a better understanding of what is contributing to these higher values, a further analysis of the cost of $\Pcf$ has to be carried out. Besides the total time spent to interpolate the coarse finite element field $\delta^{l-1}$
from $\mathcal{T}_{l-1}$ on $\mathcal{T}_{l}$, we also consider the time spent internally in evaluating $\delta^{l}$ at arbitrarily located reference points determined after the search procedure explained in Sec.~\ref{subsubsec:Search procedure}. Even if in these particular instances levels are nested and hence points are not truly arbitrary, this is transparent to the algorithmic realization
of the method and gives an indication about how much time is spent in sending messages compared to the effective evaluation cost. To measure this, we consider once again the $\mathtt{cube}$ and $\mathtt{ball}$ examples, but with $l=8$ refinement levels and polynomial degree $p=4$, in order to appreciate the different contributions. The findings are reported in Figure~\ref{fig:breakdown_evaluations}, which shows that evaluation of the polynomial functions is the dominant component of prolongation kernels, taking roughly
the $80\%$ of the total time for the present setup and geometries when using the optimized implementation presented in~\cite{Bergbauer2024}, whereas the remaining part can be associated with communication. In Figure \ref{fig:exclusive_time} we show the minimum and maximum time spent on each level for the $\mathtt{cube}$ test case for a parallel run with $12$ processes for both approaches. With a similar distribution between levels, the non-nested variant exhibits higher computational times overall because of the cost of the intergrid transfers. Notably, both variants
exhibit a reduction of times by a factor of approximately $8$ for the finest levels, as predicted by the theory.

Figure~\ref{fig:exclusive_time} also confirms that the good workbalance for the time to solution typical of the global coarsening algorithms is seamlessly inherited by the non-nested approach. We have verified that increasing the degree of the Chebyshev smoothers in order to reduce the impact of the level transfer by possibly reducing the total number of conjugate gradient iterations does not improve the overall time to solution.

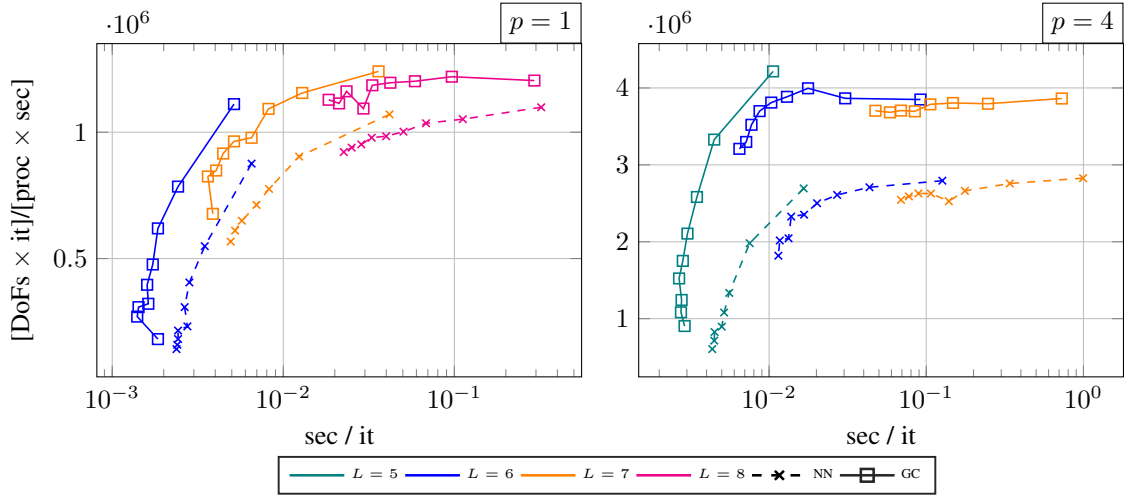
\begin{figure}[ht]
    \centering

    \vspace{.2cm}
        \begin{minipage}[t]{0.5\textwidth}
            \begin{tikzpicture}
            \begin{semilogxaxis}[
                title style={font=\tiny}, every axis title/.style={above left, at={(1,1)}, draw=black, fill=white},
                title={$p=1$},
                xlabel={sec / it},
                ylabel={[DoFs $\times$ it]/[proc $\times$ sec]},
                legend style={at={(1.05,1)}, anchor=north west},
                width=\textwidth,
                height=6cm,
                grid
            ]
        
            \addplot[solid, color=blue, line width=0.6pt, mark=square, mark options={solid, fill=none}] coordinates 
            {
                (0.00515, 1111100.0)
                (0.002425, 784770.0)
                (0.00185, 619410.0)
                (0.001725, 476630.0)
                (0.0016, 396340.0)
                (0.001625, 321210.0)
                (0.001425, 307960.0)
                (0.0014, 270410.0)
                (0.00185, 181420.0)
            };
        
            \addplot[dashed, color=blue, line width=0.6pt, mark=x, mark options={solid, fill=none}] coordinates 
            {
                (0.006525, 875750.0)
                (0.003475, 549090.0)
                (0.002825, 405480.0)
                (0.00265, 307930.0)
                (0.00275, 231680.0)
                (0.002425, 215270.0)
                (0.002425, 181730.0)
                (0.0024, 158870.0)
                (0.002375, 141890.0)
            };
        
            \addplot[solid, color=orange, line width=0.6pt, mark=square, mark options={solid, fill=none}] coordinates 
            {
                (0.03605, 1240600.0)
                (0.0129, 1155200.0)
                (0.0082, 1092300.0)
                (0.006525, 977950.0)
                (0.00515, 963680.0)
                (0.00445, 915610.0)
                (0.00405, 848360.0)
                (0.003625, 824780.0)
                (0.003875, 677010.0)
            };
        
            \addplot[dashed, color=orange, line width=0.6pt, mark=x, mark options={solid, fill=none}] coordinates 
            {
                (0.041775, 1070800.0)
                (0.012375, 903660.0)
                (0.00825, 775590.0)
                (0.006975, 712420.0)
                (0.005725, 649850.0)
                (0.005225, 610650.0)
                (0.004925, 567100.0)
            };
        
            \addplot[solid, color=magenta, line width=0.6pt, mark=square, mark options={solid, fill=none}] coordinates 
            {
                (0.293625, 1204300.0)
                (0.09665, 1219500.0)
                (0.05885, 1201600.0)
                (0.04225, 1195700.0)
                (0.03315, 1185600.0)
                (0.0294, 1093500.0)
                (0.023425, 1161700.0)
                (0.021175, 1113900.0)
                (0.018425, 1128300.0)
            };
        
            \addplot[dashed, color=magenta, line width=0.6pt, mark=x, mark options={solid, fill=none}] coordinates 
            {
                (0.3219, 1098600.0)
                (0.1121, 1051600.0)
                (0.0683, 1035400.0)
                (0.050425, 1002000.0)
                (0.039925, 984110.0)
                (0.032875, 978170.0)
                (0.028575, 951810.0)
                (0.0251, 938990.0)
                (0.022575, 921480.0)
                
            };
        
            \end{semilogxaxis}
            \end{tikzpicture}
        \end{minipage}%
        \hfill
        \begin{minipage}[t]{0.5\textwidth}
            \begin{tikzpicture}
                \begin{semilogxaxis}[
                    title style={font=\tiny},every axis title/.style={above left,at={(1,1)},draw=black,fill=white},
                    title={$p=4$},
                    xlabel={sec / it},
                    legend style={at={(1.05,1)}, anchor=north west},
                    width=\textwidth,
                    height=6cm,
                    grid
                ]

            \addplot[solid, color=teal, line width=0.6pt, mark=square, mark options={solid, fill=none}] coordinates 
            {
                (0.0106, 4215300.0)
                (0.004475, 3328300.0)
                (0.003475, 2582400.0)
                (0.003025, 2106400.0)
                (0.002825, 1751600.0)
                (0.002675, 1523200.0)
                (0.002775, 1243400.0)
                (0.00275, 1083000.0)
                (0.0029, 904520.0)

            };
                        
            \addplot[dashed, color=teal, line width=0.6pt, mark=x, mark options={solid, fill=none}] coordinates 
            {
                (0.0166, 2693700.0)
                (0.007525, 1981300.0)
                (0.005575, 1335600.0)
                (0.005175, 1080400.0)
                (0.005, 895960.0)
                (0.0045, 826890.0)
                (0.004475, 714720.0)
                (0.00435, 604670.0)
            };
                        
            \addplot[solid, color=blue, line width=0.6pt, mark=square, mark options={solid, fill=none}] coordinates 
            {
                (0.091875, 3849500.0)
                (0.0305, 3864900.0)
                (0.0177, 3995500.0)
                (0.013, 3884900.0)
                (0.010325, 3809800.0)
                (0.0087, 3698800.0)
                (0.007725, 3520900.0)
                (0.00715, 3298000.0)
                (0.006475, 3208500.0   )
            };
                        
            \addplot[dashed, color=blue, line width=0.6pt, mark=x, mark options={solid, fill=none}] coordinates 
            {
                (0.126575, 2793900.0)
                (0.043525, 2708700.0)
                (0.0271, 2609200.0)
                (0.0202, 2501300.0)
                (0.0167, 2352200.0)
                (0.013825, 2327400.0)
                (0.0133, 2045100.0)
                (0.0117, 2016700.0)
                (0.01145, 1818400.0)
            };
                        
            \addplot[solid, color=orange, line width=0.6pt, mark=square, mark options={solid, fill=none}] coordinates 
            {
                (0.728175, 3862600.0)
                (0.24705, 3794900.0)
                (0.14785, 3804800.0)
                (0.1061, 3786800.0)
                (0.08455, 3696000.0)
                (0.069, 3705600.0)
                (0.05875, 3682700.0)
                (0.047475, 3703600.0)

            };
                        
            \addplot[dashed, color=orange, line width=0.6pt, mark=x, mark options={solid, fill=none}] coordinates 
            {
                (0.9947, 2827600.0)
                (0.339975, 2757600.0)
                (0.176125, 2661700.0)
                (0.139125, 2527100.0)
                (0.10705, 2627100.0)
                (0.089175, 2628400.0)
                (0.07755, 2590300.0)
                (0.069075, 2545000.0   )
            };

        \end{semilogxaxis}
        \end{tikzpicture}
        \end{minipage}
        \hspace{0.6cm}\begin{tikzpicture}[thick,scale=1.2, every node/.style={scale=0.85}]
            \begin{axis}[%
                hide axis,
                legend style={font=\tiny},
        xmin=10,
        xmax=1000,
        ymin=0,
        ymax=0.4,
        semithick,
        legend style={draw=white!15!black,legend cell align=left},legend columns=-1
        ]
    
        \addlegendimage{teal,every mark/.append style={fill=teal!80!black}}
        \addlegendentry{$L=5$};
        \addlegendimage{blue,every mark/.append style={fill=blue!80!black}}
        \addlegendentry{$L=6$};
        \addlegendimage{orange,every mark/.append style={fill=orange!80!black}}
        \addlegendentry{$L=7$};
        \addlegendimage{magenta,every mark/.append style={fill=magenta!80!black}}
        \addlegendentry{$L=8$};
        
        \addlegendimage{dashed,black,every mark/.append style={fill=black,solid},mark=x}
        \addlegendentry{NN};
        \addlegendimage{mark=square}
        \addlegendentry{GC};

    \end{axis}
\end{tikzpicture}%
    \caption{Comparison of normalized throughput for global coarsening (GC) and non-nested (NN) algorithm applied to the $\mathtt{cube}$ test case for polynomial degrees $p=1$ (left) and $p=4$ (right).}
    \label{fig:throughput_GC_NN_cube}
\end{figure}

\noindent Finally, Figure~\ref{fig:throughput_GC_NN_cube} displays the results of a strong scaling test starting from 1 compute node ($48$ processes) up to 16 compute nodes ($768$ processes), applied to the $\mathtt{cube}$ test case with polynomial degrees $p=1$ and $p=4$ and different refinement levels $l$. In particular, we plot the \emph{normalized throughput} (computed as the number of degrees of freedom per process and time per iteration) against the \emph{time per iteration} both for the global coarsening and non-nested algorithms. As expected, the throughput achieved
for $p=4$ is higher than for $p=1$. This does not just hold for global coarsening (as already reported in~\cite{MHPSK}), but also for the non-nested methodology, thanks to the usage of matrix-free algorithms within smoothers, which are known to improve in throughput with the polynomial order~\cite{Kronbichler2018}. The gap between curves associated to the same refinement level $l$ is more noticeable in favor of global coarsening with $p=4$.

\begin{figure}[ht]
    \centering
    \resizebox{\columnwidth}{!}{%
    \begin{tikzpicture}
        \pgfplotsset{width=12cm}
        \begin{axis}
            [
            title=\Large{Cube},
            ymin=0,
            ymax=.8,
            bar width=0.15,
            legend style={at={(0.4,-0.25)}, 
            anchor=north,legend columns=-1},
            xlabel={Level},
            ylabel={Time [s]},
            xtick={6,7,8},
            xticklabels={6,7,8},
            grid=major,
            grid style={line width=.1pt, draw=gray!50},
            legend pos=north west,
            ]

            \addlegendimage{ybar,ybar legend,blue, fill=blue!40!white}
            \addlegendentry{Evaluation};
            \addlegendimage{ybar,ybar legend,red, fill=red!40!white}
            \addlegendentry{Total prolongation};
            \addplot [ybar stacked, blue, fill=blue!40!white] coordinates {
            (6, 6.4654e-03)
            (7, 5.1699e-02)
            (8, 4.1978e-01)
            };
            \addplot [ybar stacked, red, fill=red!40!white] coordinates {
            (6, 9.2360e-03 - 6.4654e-03)
            (7, 6.9205e-02 - 5.1699e-02)
            (8, 5.5324e-01 - 4.1978e-01)
            };

        \end{axis} 
    \end{tikzpicture}
    \begin{tikzpicture}
        \pgfplotsset{width=12cm}
        \begin{axis}
            [
            title= \Large{Ball},
            ymin=0,
            ymax=.8,
            bar width=0.15,
            legend style={at={(0.4,-0.25)}, 
            anchor=north,legend columns=-1},
            xlabel={Level},
            xtick={6,7,8},
            xticklabels={6,7,8},
            grid=major,
            grid style={line width=.1pt, draw=gray!50},
            legend pos=north west,
            ]

            \addlegendimage{ybar,ybar legend,blue, fill=blue!40!white}
            \addlegendentry{Evaluation};
            \addlegendimage{ybar,ybar legend,red, fill=red!40!white}
            \addlegendentry{Total prolongation};
            \addplot [ybar stacked, blue, fill=blue!40!white] coordinates {
            (6, 4.5237e-02)
            (7, 5.1763e-02)
            (8, 4.1850e-01)
            };
            \addplot [ybar stacked, red, fill=red!40!white] coordinates {
            (6, 6.2375e-02 - 4.5237e-02)
            (7, 6.9496e-02 - 5.1763e-02)
            (8, 5.5341e-01 - 4.1850e-01)
            };
        \end{axis} 
    \end{tikzpicture}
    }
    \caption{Breakdown of the computational time (measured in seconds) required for prolongation between levels with $12$ processes. In red: total time spent to prolongate a coarse finite element field from level $l$ to $l+1$. In blue: time spent entirely on the evaluation at arbitrarily located reference points. For both the $\mathtt{cube}$ (left) and $\mathtt{ball}$ (right) meshes, we perform the experiment with polynomial degree $p=4$ and $l=8$ refinement levels}
    \label{fig:breakdown_evaluations}
\end{figure}
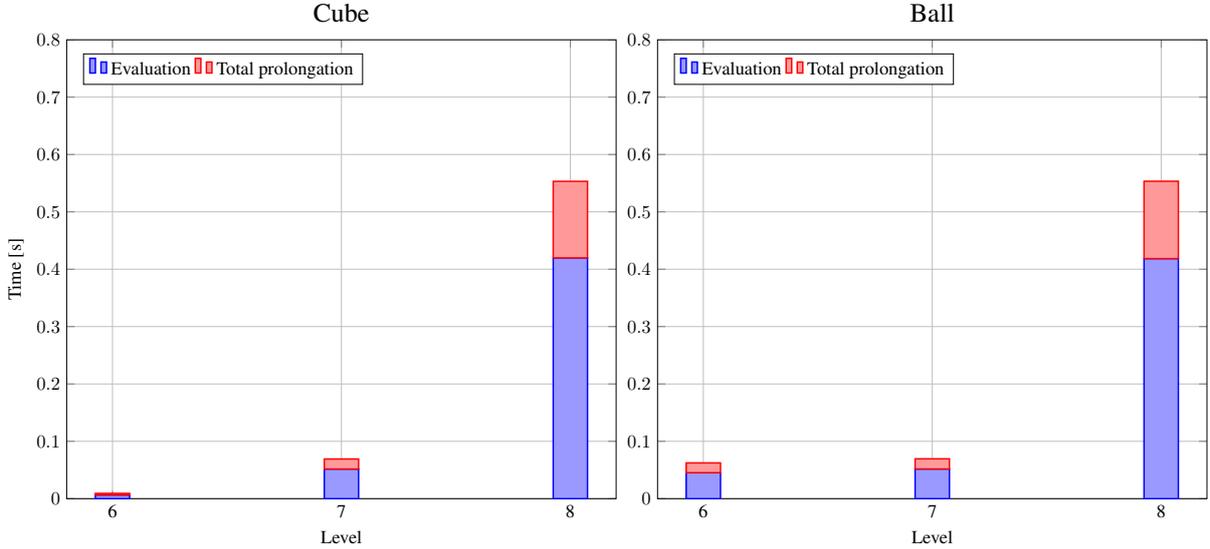

\subsection{Non-nested meshes: performance metrics and mesh partitioning}
In this subsection, we analyze the two non-nested tests presented in Section \ref{subsub:L_shaped} more in depth. An analysis of the different components required by the non-nested multigrid method is presented. We adapt the setup and tests described in \cite{MHPSK,CHKK}. Almost every metric adopted therein
and in classical multigrid literature can be immediately reused, except for the \emph{vertical communication efficiency}, i.e., the share of fine cells that are owned by the same process that owns their parent coarse cell. That metric quantifies the amount
of data (such as number of points) that has to be exchanged by during the level transfer. The lack of exact overlapping between distributed levels slightly complicates standard metrics that are derived solely on
geometrical information, usually considered in standard multigrid literature. In order to extend this definition to our context, we replace the concept of children cells with the one
of \emph{owned points}, so that the definition becomes:
\begin{definition}[Vertical communication efficiency for non-matching levels]
    \label{def:vertical_communication}
    Given two arbitrarily partitioned triangulations discretizing $\Omega$, the \emph{vertical communication efficiency} is the share of \emph{owned points} on the finer grid that have the same
    owning process as their corresponding owners on the coarse grid.
\end{definition}
\hfill \break
We consider the following geometrical metrics, taken from~\cite{MHPSK,CHKK}:
\begin{itemize}
    \item \textbf{Serial workload}: Sum of the number of cells on all levels $\mathcal{W}_s \coloneqq \sum_{l} \mathcal{C}_l$, being $\mathcal{C}_l$ the number of cells on the $l$-th level.
    \item \textbf{Parallel workload}: Sum of the maximum number of cells owned by any process on every level of the hierarchy: $\mathcal{W}_p \coloneqq \sum_l \max_{p}{ \mathcal{C}_{l}^{p}}$, where $\mathcal{C}_{l}^{p}$ is
    the number of cells owned by process $p$ on level $l$. Load imbalances imply that in general one has $\mathcal{W}_p \ne \frac{\mathcal{W}_s}{p}$, meaning that the work is not properly distributed on the levels. The \textbf{parallel workload efficiency} is hence defined as $\frac{\mathcal{W}_s}{\mathcal{W}_p \cdot p}$.
    \item \textbf{Vertical efficiency} according to Definition~\ref{def:vertical_communication}: This metric gives a good indication on how much data has to be exchanged in order to perform intergrid transfers. A large value means that consecutive partitions of the
          levels are well overlapped, implying a small volume of communication. Within this work, all meshes are partitioned independently by employing the default partitioning strategy of~\textsc{deal.II}. Optionally, graph partitioners such as $\mathtt{METIS}$~\cite{METIS} can be employed for such a task.
\end{itemize}

\begingroup
\renewcommand{\arraystretch}{1.1}
\begin{table}[!t]
    \centering
    \begin{minipage}[t]{0.45\textwidth}
        \centering
        \begin{tabular}{l|ccc}
            \toprule
            \multicolumn{4}{c}{$12$ processes (default partitioning policy)} \\
            \cmidrule(lr){1-4}
            \emph{l} & wl & wl-eff & v-eff \\
            \midrule
            3 & \num{8.5}e+2  & 99\% & 7\%\\
            4 & \num{1.5}e+3  & 99\% & 6\% \\
            5 & \num{2.7}e+3  & 99\% & 15\% \\
            6 & \num{3.9}e+3  & 99\% & 39\%\\
            \bottomrule
        \end{tabular}
        \\[1ex]
        \makebox[0.95\textwidth][c]{(a) Default partitioning policy.}
    \end{minipage}
    \hspace{0.5cm}
    \begin{minipage}[t]{0.45\textwidth}
        \centering
        \begin{tabular}{l|ccc}
            \toprule
            \multicolumn{4}{c}{$12$ processes  (matching partitioning policy)} \\
            \cmidrule(lr){1-4}
            \emph{l} & wl & wl-eff & v-eff \\
            \midrule
            3 & \num{1.0}e+3 & 85\% & 89\%\\
            4 & \num{1.9}e+3 & 79\% & 91\% \\
            5 & \num{3.3}e+3 & 81\% & 92\% \\
            6 & \num{4.5}e+3 & 86\% & 92\%\\
            \bottomrule
        \end{tabular}
        \\[1ex]
        \makebox[0.95\textwidth][c]{(b) Matching partitioning policy.}
    \end{minipage}
    \\[2ex]
    \caption{Multigrid statistics for the 2D L-shaped test case for different number of levels (wl: parallel workload, wl-eff: parallel workload efficiency, v-eff: vertical communication efficiency.) (a) Default partitioning policy. (b) Matching partitioning policy.}
    \label{tab:lshape_metrics}
\end{table}
\endgroup

\begingroup
\renewcommand{\arraystretch}{1.1}
\begin{table}[!t]
    \centering
    \begin{minipage}[t]{0.45\textwidth}
        \centering
        \begin{tabular}{l|ccc}
            \toprule
            \multicolumn{4}{c}{$12$ processes (default partitioning policy)} \\
            \cmidrule(lr){1-4}
            \emph{l} & wl & wl-eff & v-eff \\
            \midrule
            2 & \num{3.7}e+2 & 99\% & 1\% \\
            3 & \num{1.2}e+3 & 99\% & 2\% \\
            4 & \num{3.5}e+3 & 99\% & 5\% \\
            \bottomrule
        \end{tabular}
        \\[1ex]
        \makebox[0.95\textwidth][c]{(a) Default partitioning policy.}
    \end{minipage}
    \hspace{0.5cm}
    \begin{minipage}[t]{0.45\textwidth}
        \centering
        \begin{tabular}{l|ccc}
            \toprule
            \multicolumn{4}{c}{$12$ processes  (matching partitioning policy)} \\
            \cmidrule(lr){1-4}
            \emph{l} & wl & wl-eff & v-eff \\
            \midrule
            2 & \num{5.1}e+2 & 72\% & 60\% \\
            3 & \num{2.2}e+3 & 58\% & 67\% \\
            4 & \num{6.6}e+3 & 53\% & 72\% \\
            \bottomrule
        \end{tabular}
        \\[1ex]
        \makebox[0.95\textwidth][c]{(b) Matching partitioning policy.}
    \end{minipage}
    \\[2ex]
    \caption{Multigrid statistics for the 3D Fichera test case for different number of levels (wl: parallel workload, wl-eff: parallel workload efficiency, v-eff: vertical communication efficiency.) (a) Default partitioning policy. (b) Matching partitioning policy.}
    \label{tab:Fichera_metrics}
\end{table}
\endgroup

\begingroup
\renewcommand{\arraystretch}{1.1}
\begin{table}[!t]
    \centering
    \begin{minipage}[t]{0.45\textwidth}
        \centering
        \begin{tabular}{l|ccc}
            \toprule
            \multicolumn{4}{c}{$12$ processes (default partitioning policy)} \\
            \cmidrule(lr){1-4}
            \emph{l} & wl & wl-eff & v-eff \\
            \midrule
            2 & \num{1.5}e+3 & 99\% & 1\%\\
            3 & \num{4.1}e+3 & 99\% & 2\%\\
            4 & \num{7.9}e+3 & 99\% & 5\%\\
            \bottomrule
        \end{tabular}
        \\[1ex]
        \makebox[0.95\textwidth][c]{(a) Default partitioning policy.}
    \end{minipage}
    \hspace{0.5cm}
    \begin{minipage}[t]{0.45\textwidth}
        \centering
        \begin{tabular}{l|ccc}
            \toprule
            \multicolumn{4}{c}{$12$ processes  (matching partitioning policy)} \\
            \cmidrule(lr){1-4}
            \emph{l} & wl & wl-eff & v-eff \\
            \midrule
            2 & \num{1.6}e+3 & 96\% & 87\%\\
            3 & \num{4.4}e+3 & 92\% & 88\%\\
            4 & \num{8.3}e+3 & 95\% & 89\%\\
            \bottomrule
        \end{tabular}
        \\[1ex]
        \makebox[0.95\textwidth][c]{(b) Matching partitioning policy.}
    \end{minipage}
    \\[2ex]
    \caption{Multigrid statistics for the piston test case for different number of levels (wl: parallel workload, wl-eff: parallel workload efficiency, v-eff: vertical communication efficiency.) (a) Default partitioning policy. (b) Matching partitioning policy.}
    \label{tab:piston_metrics}
\end{table}
\endgroup

\begin{figure}[h!]
    \centering
    \begin{minipage}[t]{0.45\textwidth}
        \centering
        \includegraphics[width=\textwidth]{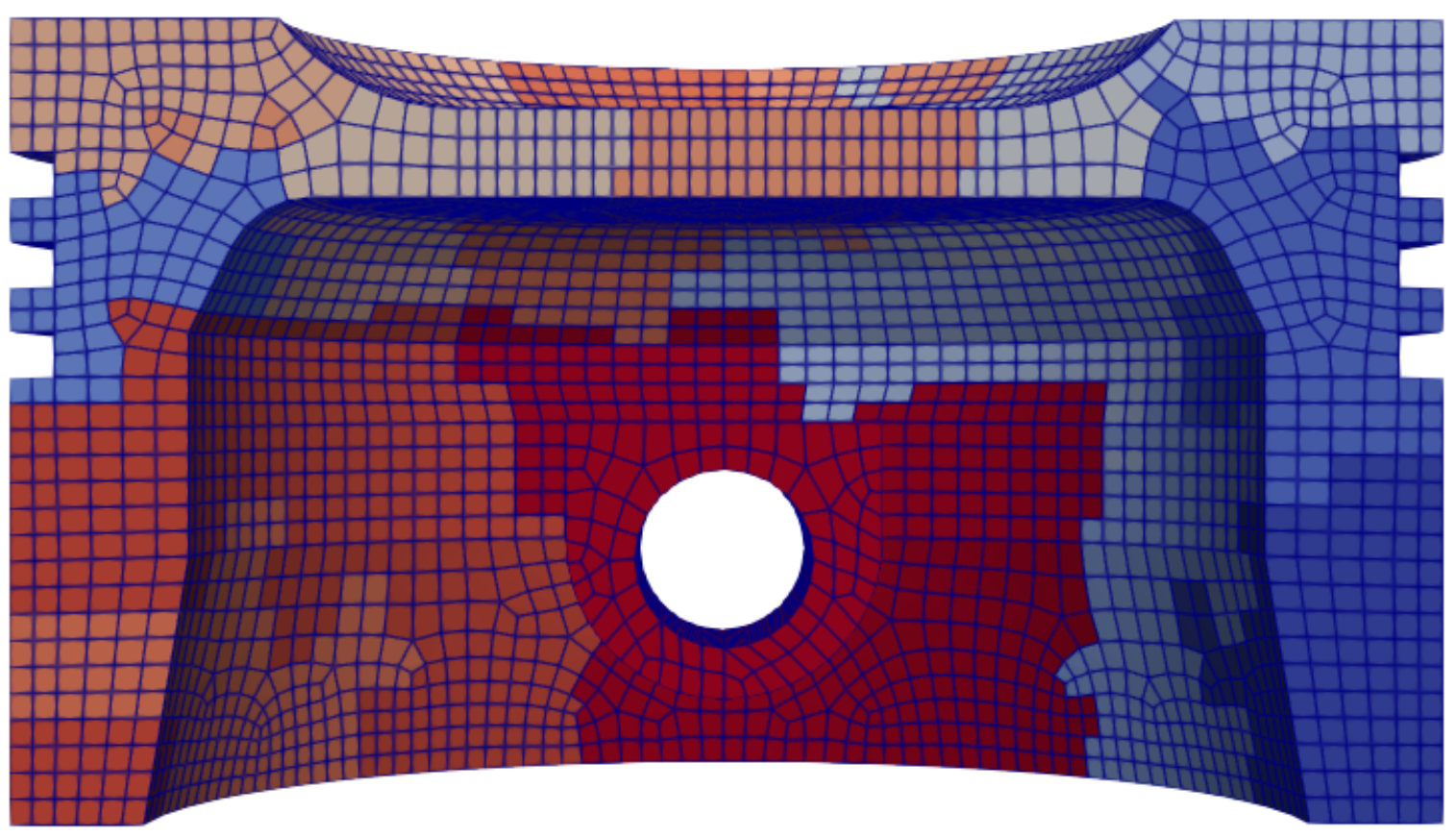}
        \caption{Partitioning for fine level.}
        \label{fig:coarse_piston_original}
    \end{minipage}
    \hfill
    \begin{minipage}[t]{0.45\textwidth}
        \centering
        \includegraphics[width=\textwidth]{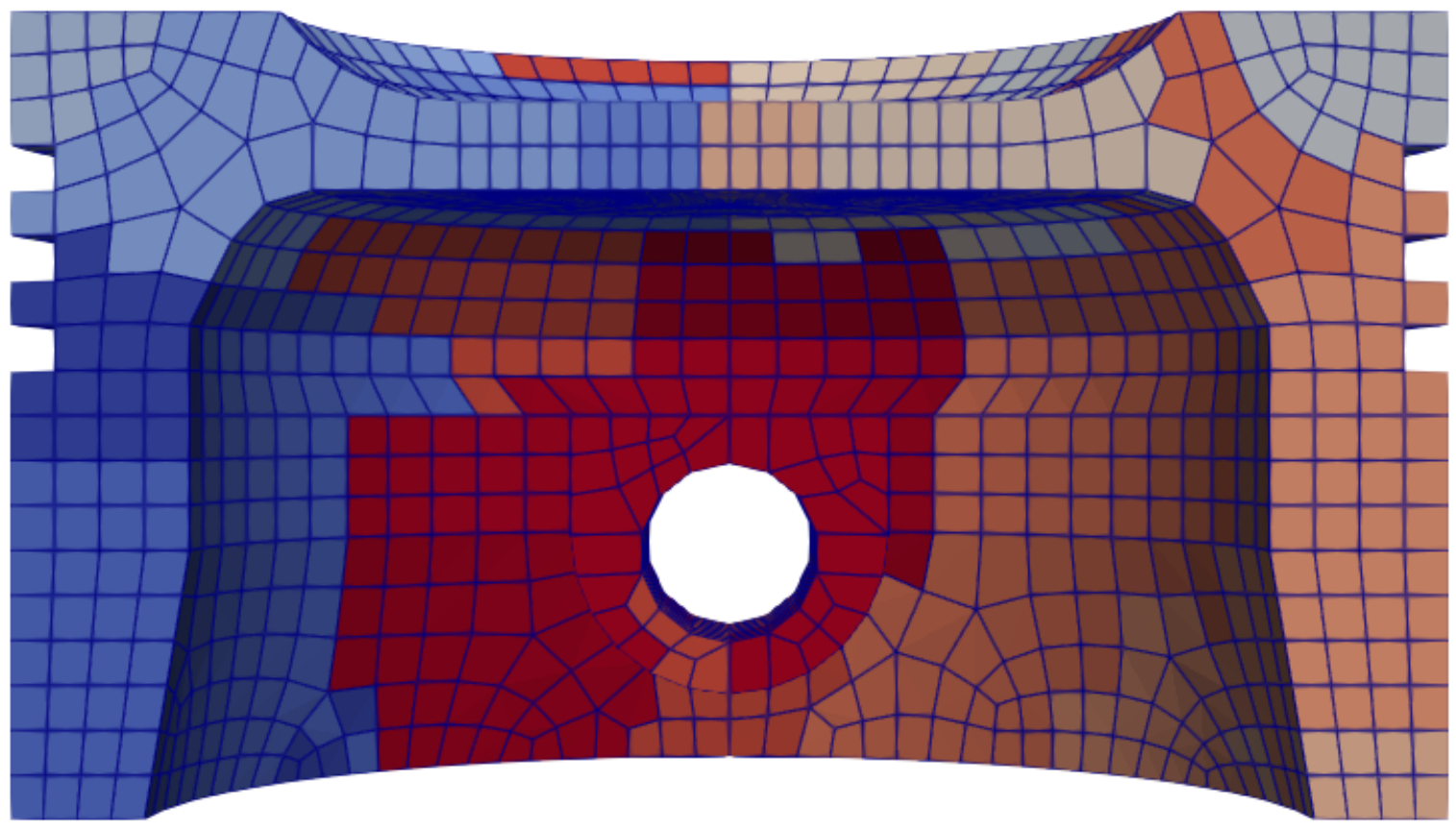}
        \caption{Coarse level (default partitioning).}
        \label{fig:fine_piston}
    \end{minipage}
    
    \vspace{0.1cm} %
    \begin{minipage}[t]{0.5\textwidth}
        \centering
        \includegraphics[width=\textwidth]{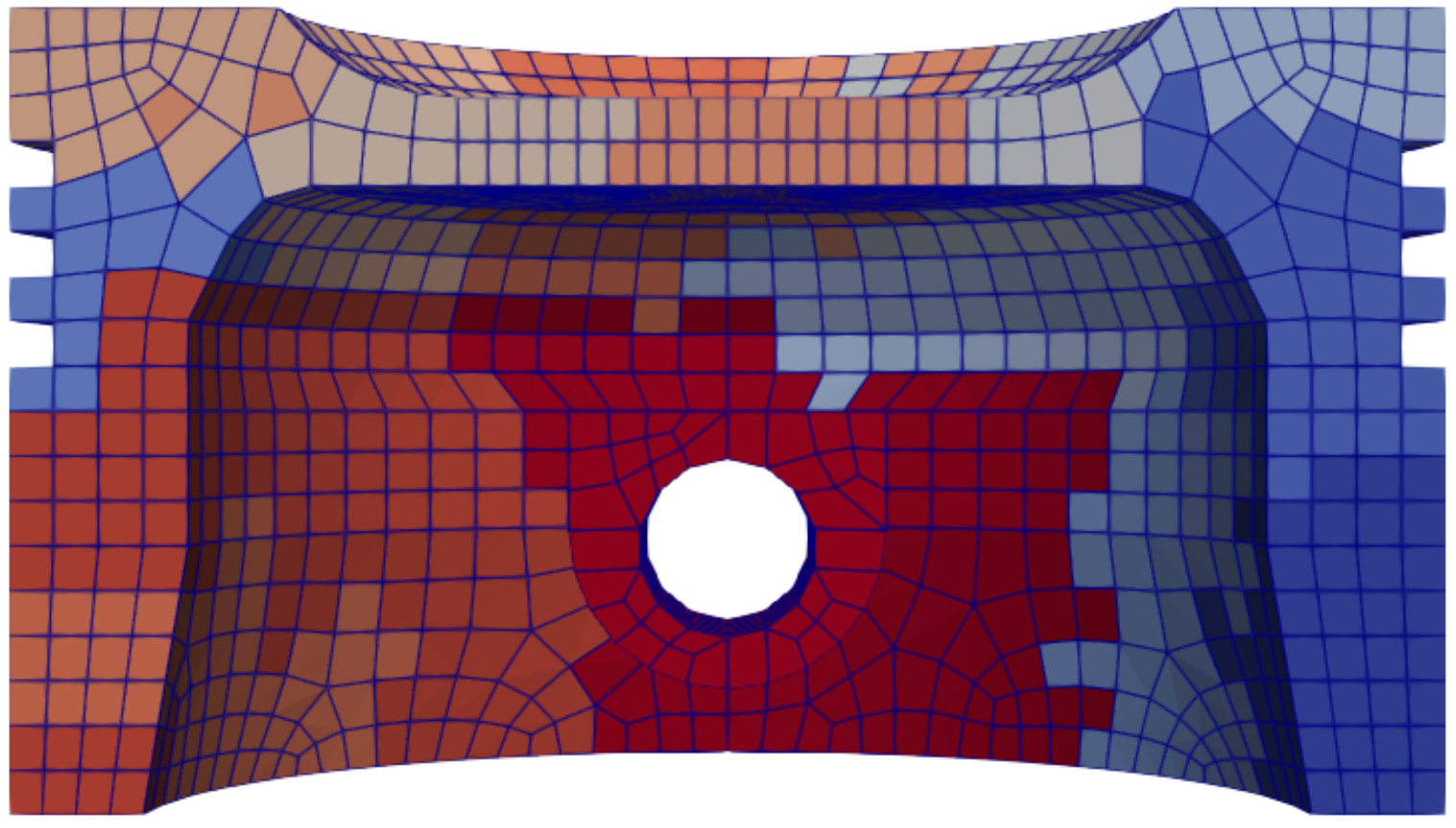}
        \caption{Coarse level (repartitioned).}
        \label{fig:coarse_piston_repatitioned}
    \end{minipage}
    
    \caption{Top row: clipped view of consecutive levels in the hierarchy for the piston test case~\ref{fig:piston} partitioned across $12$ processors without any partioning policy (each color represents a different MPI rank).
    Bottom row: repartitioning of the coarse level. Notice the high match in terms of colors with the fine level on the top-left and the different mesh-sizes.}
    \label{fig:processors_coupling_piston}
\end{figure}

\noindent In the global coarsening case the parallel workload is, by construction, well distributed among the participating processes as each level is partitioned during the construction of the hierarchy. On the other hand, this affects the value of vertical efficiency which
is generally low. Indeed, an optimal load balance and an optimal overlap of parallel partitions are mutually orthogonal requirements. Overall, this implies a higher computational pressure on the data exchange phase needed for prolongation and restriction operations. In this context, where levels do not share the same coarse mesh, the drop in vertical efficiency
is even more pronounced compared to global coarsening. Tables~\ref{tab:lshape_metrics}(a),~\ref{tab:Fichera_metrics}(a) confirm this expected behavior for the levels displayed in Figures \ref{fig:levels_L_shape} and \ref{fig:levels_Fichera}. In particular, the good parallel workload efficiency of global coarsening
is automatically inherited, at the cost of a poor vertical efficiency for both geometries. Table \ref{tab:piston_metrics}(a) reports the same statistics for the piston test case in \ref{fig:piston}, partitioned with $12$ processors. The statements made for
the previous geometries still hold also for the piston hierarchy. In this case the low value of the vertical efficiency is even more prominent. 
As it can be appreciated from the top row in Figure~\ref{fig:processors_coupling_piston}, the overlap of the parallel partitions for two consecutive levels
is quite low, and
 our findings match with patterns reported in \cite{MHPSK}, so making definite statements is not possible as they would depend on the number of processes and are in practice problem-specific.
A possible workaround to increase the vertical efficiency between levels could consist in a different partitioning policy which aims to match the partitioning on the finest (and computationally most expensive) levels of the hierarchy, instead of relying on completely independently distributed hierarchies. We investigate this
by first extracting the partitioning of the finest triangulation of the hierarchy, and then repartitioning coarser meshes after having determined the processes on the fine mesh owning portions of each coarse mesh. The outcome of such approach is depicted in Figure~\ref{fig:coarse_piston_repatitioned} for the piston test case, while
Tables~\ref{tab:lshape_metrics}(b), \ref{tab:Fichera_metrics}(b), and \ref{tab:piston_metrics}(b) illustrate the result of the application of this \emph{matching} policy.  We observe a large increase in the vertical efficiency, without dramatic effects on the parallel workload. Notably, in the piston test case we achieve an average vertical efficiency of $88\%$ while keeping nearly optimal
workload efficiency. The advantage of the induced higher vertical efficiency would be a decrease of the setup cost associated with the geometric search procedure outlined in Section~\ref{sec:implementation}.

